\documentclass[12pt]{article}

\usepackage{amsfonts,amsmath,amsthm,latexsym,amssymb,fullpage,amscd}
\usepackage[mathscr]{eucal}

\newcommand{\R}{\mathbb{R}}

\newtheorem{theorem}{Theorem}[section]
\newtheorem{definition}[theorem]{Definition}
\newtheorem{lemma}[theorem]{Lemma}
\newtheorem{corollary}[theorem]{Corollary}
\title{Smooth and singular maximal averages over 2D hypersurfaces and associated Radon transforms}
\author{Michael Greenblatt}
\date{\today}
     
\newcommand\blfootnote[1]{%
  \begingroup
  \renewcommand\thefootnote{}\footnote{#1}%
  \addtocounter{footnote}{-1}%
  \endgroup
}
\begin{document}
\maketitle
\begin{abstract}  We prove $L^p$ boundedness results, $p > 2$, for local maximal averaging operators over a smooth 2D hypersurface $S$ with 
either a $C^1$ density function or a density function with a singularity that grows as $|(x,y)|^{-\beta}$ for $\beta < 2$. 
Suppose one is in coordinates such that the 
surface is localized near some $(x_0,y_0,z_0)$ at which $(0,0,1)$ is normal to the surface, and suppose the surface is
represented as the graph
of $z_0 + s(x - x_0, y - y_0)$ near $(x_0,y_0)$, with $s(0,0) = 0$. It is shown that as long as the Taylor series of the
Hessian determinant of $s(x,y)$ at $(0,0)$ is not identically zero, the maximal
averaging operator is bounded on $L^p$ for $p > \max(2,1/g)$, where $g$ is an index based on the growth rate of
the distribution function $s(x,y)$ near the origin. Standard examples show that the exponent $1/g$ is best
possible whenever the tangent plane to $S$ at $(x_0,y_0,z_0)$ does not contain the origin. This theorem improves on the main result of [IKeM],
using different methods.
We use closely related methods to prove $L^p$ to $L^p_{\alpha}$ Sobolev estimates for Radon transform operators with the
same density functions, with no excluded cases. In the $g < 1/2$ case, there is an interval $I$ containing
$2$ for which $L^p$ to $L^p_{\alpha}$ boundedness is proven for $\alpha < g$ when $p \in I$, and for such $p$ one can never 
gain more than $g$ derivatives.
\end{abstract}

\section{Background and Theorem Statements}

\subsection {Maximal averaging operators}

Let $S$ be a real analytic hypersurface in $\R^3$ and let $X_0 = (x_0,y_0,z_0)$ be a point on $S$. Letting $X = (x,y,z)$ denote a point in
$\R^3$, for a real-valued function $\psi(X)$ on $\R^3$ supported near $X_0$ we will be considering two operators. The first is the maximal average
$$Mf(X) = \sup_{t > 0} \bigg|\int  f(X - tX')\, \psi(X')\, d\lambda(X') \bigg|\eqno (1.1)$$
Here $d\lambda(X')$ denotes the standard Euclidean surface measure on $S$ in the $X'$ variable. We will be examining the $L^p$ boundedness of
$M$ on $L^p$ for $p > 2$ when $\psi$ is $C^1$ or has a singularity at $(x_0,y_0,z_0)$.

\noindent The other operator we will be looking at is the Radon transform operator
$$Rf(X) = \int f(X - X')\, \psi(X')\, d\lambda(X') \eqno(1.2)$$
For $R$ we will be proving $L^p$ to $L^p_{\alpha}$ Sobolev space estimates.

 By the rotation-invariance of the above estimates, it suffices to consider the case when 
the vector $(0,0,1)$ is normal to $S$ at $(x_0,y_0,z_0)$, and we always make this assumption in this paper. Now the surface $S$ is the graph of 
$z_0 + s(x - x_0, y - y_0)$ for some real analytic $s(x,y)$ defined near $(0,0)$ satisfying
$$ s(0,0) = 0 {\hskip 0.7 in} \nabla s(0,0) = (0,0)\eqno (1.3)$$
In the new coordinates one can replace $\psi(x,y,z)$ by a function $\phi(x,y)$ of $x$ and $y$ only. The type of theorems we are proving are readily
shown to be false if $s(x,y)$ is identically zero, so we will always assume that this is not the case.
The conditions we will be assuming on $\phi(x,y)$ are as follows. Let ${\bf x}$ denote
$(x,y)$ and ${\bf x_0}$ denote $(x_0,y_0)$. We assume that $\phi(x,y)$ is $C^1$ on $\R^2 -\{{\bf x_0}\}$ and that for some $A \geq 0$ and 
$0 \leq \beta < 2$ we have
$$|\phi({\bf x})| \leq A |{\bf x} - {\bf x_0}|^{-\beta} {\hskip 0.7 in} |\nabla \phi({\bf x})| \leq A |{\bf x} - {\bf x_0}|^{-\beta - 1} \eqno (1.4)$$
The case when $\beta = 0$ includes the case when $\phi(x,y)$ is $C^1$, and in fact when $\beta =  0$ the sharp cases of our 
theorems will always correspond to the sharp estimates for the $C^1$ situation.

\noindent Using $s(x,y)$ and $\phi(x,y)$, the Radon transform operator can be rewritten in the form
$$Rf(X) = \int_{\R^2} f(x - x',y - y', z - z_0 - s(x' - x_0, y' - y_0))\, \phi(x',y')\,dx' dy' \eqno (1.5) $$
Note that by the translation-invariance of $L^p$ Sobolev estimates for Radon transforms, we can always assume that $(x_0,y_0,z_0) = (0,0,0)$ in the 
analysis of $Rf$, but this is not immediately the case with the maximal averages. We will see however that our arguments are essentially independent
of $(x_0,y_0,z_0)$ for both operators, and in fact our analysis for the two operators will be quite similar. Nonetheless, for maximal averages there 
are certain situations where the best possible result does depend on $(x_0,y_0,z_0)$. We refer to the discussion after Theorem 1.1 for more on this.

Our theorem statements will be in terms of the quantity $g > 0$ such that there exists
$d = 0$ or $1$ such that for any sufficiently small $r > 0$, as $\epsilon \rightarrow 0$ for some $C_{\beta,r} > 0$ one has asymptotics
$$\int_{\{({\bf x},z) \in S:\,|{\bf x} - {\bf x_0}| < r,\,|z - z_0| < \epsilon \}} |{\bf x} - {\bf x_0}|^{-\beta} \,d{\bf x}= C_{\beta,r}\, \epsilon^{g}
|\ln \epsilon|^{d} + o(\epsilon^{g}
|\ln \epsilon|^{d}) \eqno (1.6)$$
That such a $g$ exists follows from resolution of singularities: in terms of the function $s(x,y)$ one may write
$$\int_{\{({\bf x},z) \in S:\,|{\bf x} - {\bf x_0}| < r,\,|z - z_0| < \epsilon \}} |{\bf x} - {\bf x_0}|^{-\beta} \,d{\bf x} =
\int_{\{{\bf x}:\,|{\bf x} - {\bf x_0}| < r,\,|s({\bf x_0} - {\bf x})| < \epsilon\}} |{\bf x} - {\bf x_0}|^{-\beta}\,d{\bf x}$$
Shifting variables from ${\bf x}$ to ${\bf x} - {\bf x_0}$, this becomes
$$ \int_{\{{\bf x}:\,|{\bf x}| < r,\,|s({\bf x})| < \epsilon\}} |{\bf x}|^{-\beta} \,d{\bf x} \eqno (1.7)$$
The theory of resolution of singularities says, roughly speaking, that if $r$ is sufficiently small, then the set $\{{\bf x}: |{\bf x}| < r\}$ can be written as the 
union of finitely many subsets on each of which a sequence of coordinate changes can be performed after which $s(\bf x)$ and $|\bf x|$ are both 
effectively monomials. The Jacobian of each sequence of coordinate changes is also effectively a monomial. As a result, the integral $(1.7)$ can be
written as a sum of several terms which are effectively of the form $(1.7)$ where in place of $\bf x$ and $s(\bf x)$ we have two monomials. There is 
also a Jacobian factor, also effectively a monomial. For such a term proving an asymptotic expansion of the form $(1.6)$ is relatively straightforward, and
 one then can add the results over several terms to get an analogous asymptotic expansion for $(1.7)$. We refer to chapter 7 of [AGV] for information on  more on such matters.  

It is not hard to show that $g$ is supremum of all $t$ for which 
$|s(x,y)|^{-t}|(x,y)|^{-\beta}$ is integrable on a neighborhood of the origin, which can serve as an alternate definition for $g$.
 Also note that $g$ is maximized when
$\beta = 0$ and $s(x,y) = x^2 + y^2$, in which case $g = 1$. Hence always $g \leq 1$. 

It is worth pointing out that in the case where $\beta = 0$ and $S$ is a convex surface of 
finite line type, $(1.7)$ illustrates that the definition of $g$ generalizes the definition of an analogous index in [BrNW]. 

\noindent Our theorem regarding the maximal averages for real analytic surfaces is as follows.

\begin{theorem} 
\

\noindent {\bf a)} Suppose the Hessian determinant of $s(x,y)$ is not identically zero. Then there is a neighborhood $N$ of $(x_0,y_0)$ such that if 
$\phi(x,y)$ is supported in $N$ and $(1.3)$-$(1.4)$ are satisfied, then $M$ is bounded on $L^p$ for $p > \max({1 \over g}, 2)$.

\noindent {\bf b)} If the tangent plane to $S$ at $(x_0,y_0,z_0)$ does not contain the origin and if on some neighborhood of $(x_0,y_0)$ one has $|\phi(x,y)| > C|{\bf x} - {\bf x_0}|^{-\beta}$ for some $C > 0$, then $M$ is not bounded on $L^p$ for any 
$1 \leq p \leq {1 \over g}$.
\end{theorem}

Note that part b) of Theorem 1.1 shows that the exponent of part a) is sharp whenever $g \leq {1 \over 2}$, if the tangent plane to $S$ at $(x_0,y_0,z_0)$ does not contain the origin. Also note that part b) does not require the Hessian condition and also holds for  $g > {1 \over 2}$. In the case where 
$\phi(x,y)$ is smooth and the tangent plane to $S$ at $(x_0,y_0,z_0)$ does not contain the origin, Theorem 1.1 follows from the Acta paper
[IKeM]. When this tangent plane does contain the origin, there are situations for which one gets stronger results. For example, in [Z] it is
 shown that
for a smooth density and $(x_0,y_0,z_0) = (0,0,0)$, $M$ is bounded on $L^p$ for all $p > 2$. On the other hand, in [Z] it is also shown that there
are some smooth surfaces containing the origin where the exponent of Theorem 1.1a) is best possible.

Why one might be interested in situations where the tangent plane to $S$ at $(x_0,y_0,z_0)$ contains the origin can be seen as follows. Suppose
$S$ is a compact surface such that the origin is contained in the exterior of $S$. Then there will typically be a one-dimensional subset $K$ of $S$ such
the tangent plane to $S$ at each point in $K$ contains the origin. Thus a result such as [IKeM] not covering such cases cannot be used to give a
boundedness result for such an $S$. When the surface everywhere has at least one nonvanishing principal curvature one can use a theorem such as [So], 
but when there are points where both principal curvatures vanish, both curvatures may vanish on a subset of 
$S$ which intersects $K$. Thus an alternative to a result such as [IKeM] is needed here.

It should be pointed out that the author wrote an earlier paper [G7] on maximal averages for two-dimensional hypersurfaces with smooth density functions.
Like this paper, the paper [G7] required that the Hessian determinant of $s(x,y)$ not be identically zero. But in addition, a second larger class of surfaces was excluded, which included examples such as $s(x,y) = (y + x^a)^b + x^c$ for $c \geq b > 2$, $a > 1$, and $ab < c$.
This class of surfaces was defined in a rather technical way using adapted coordinates and Newton polygons of certain functions derived from $s(x,y)$.
Adapted coordinate systems are coordinate systems
in which the indices $o$ and $e$ can be read off in a natural way in terms of the Newton polygon of $s(x,y)$ at the origin, as was first observed in [V].
These coordinate systems were used in an essential way in both [IKeM] and [G7]. 
The analysis for this omitted class of surface is especially difficult, and the methods the author used in [G7] involving adapted coordinates, Newton distances, 
and so on, were not easily amenable to this type of 
surface. 

In this paper, we use no such coordinate systems. Instead, unlike [IKeM] and [G7] we use full-fledged resolution of singularities theorems, in the forms
given in section 2, in conjunction with technical lemmas such as Lemma 2.5 and 2.6. As a result, we are able to go beyond what can be obtained by
 the methods of [G7] and also do not require a transversality condition as the paper [IKeM] does. We only omit certain surfaces whose analysis requires methods connected to the proof of the circular maximal theorem (many of which can be dealt with directly using such methods). This will be explained in more detail below. The robustness of our methods
is illustrated in the fact that they immediately extend to the case of 
singular density functions, where the correct notions of adapted coordinates, height, and so on would have to at least be reformulated for the singular case
when $\beta$ is large.

Another key difference between our methods and those of [IKeM] is our use of damping functions in conjunction with interpolation using a lemma 
from [SoS] (Theorem 3.1 of this paper). This enables us to avoid introducing square functions that would add significantly to the technical complexity of 
this paper, and instead reduces much of our effort to estimating oscillatory integrals.

\noindent {\bf The Hessian Condition.}

In Lemma 3.4 we will see that if the Hessian determinant of a smooth $s(x,y)$ vanishes to infinite order at the origin, then there is an invertible linear map 
$L$ such that the Taylor series of $s \circ L(x,y)$ at the origin is of the form
the form $a(x,y)y^n$ where $a(0,0) \neq 0$. In  fact (Corollary 2.2 of [dBvdE]) if $s(x,y)$ is a polynomial then $s \circ L(x,y)$ can even be written in the form $a(y)y^n$ with $a(0) \neq 0$. As a result, in the case where the Hessian vanishes to infinite order at the origin the analysis of $M$ becomes very related to the analysis of maximal averages over curves in $\R^2$. Note that these exceptional situations never occur when the surface $S$ is a compact 
surface whose defining functions are real analytic; in these exceptional situations the intersection of $S$ with the plane $y = y_0$ contains a line segment parallel to the $x$-axis, which for such a surface can only happen when this intersection contains the entire line containing this segment, which is 
not possible for a compact surface.

On the other hand, the proofs of this paper
 do not use the methods normally used to deal with maximal averages over curves in the plane. So in a sense
Theorem 1.1 a) covers all cases except those which are closely connected to maximal averages over curves in the plane. Fortunately, many of the remaining 
cases can be proved directly using existing theorems on maximal averages over curves in the plane. For example, suppose $S$ is a smooth surface 
of one of the above exceptional forms, $\phi(x,y)$ is smooth, and the
tangent plane to $S$ at $(x_0,y_0,z_0)$ does not contain the origin. Since we have rotated coordinates so that $(0,0,1)$ is normal to $S$ at $(x_0,y_0,z_0)$,
we have that
$z_0 \neq 0$. Then immediately from the definitons, one has $ Mf(X) \leq C\int_{\theta_0}^{\theta_1} M_{\theta}f(X) \,d\theta$, where $M_{\theta}f(X)$ is as
follows. Let $P_{\theta}$ denote the plane containing $X$ that is parallel to the $y$ axis and 
making an angle $\theta$ with the $z$-axis. Then $M_{\theta}f(X)$ denotes the supremum of 
the absolute values of the averages of $f$ over dilations of the curve $S \cap P_{\theta}$  that are contained in $P_{\theta}$, using the density function derived from $\phi$. In 
other words, $M_{\theta}f(X)$ is a maximal average of $f$ over curves in the 2D plane containing $X$ parallel to $P_{\theta}$.

 The $L^p$ boundedness of a given $M_{\theta}$ for $p > n$ follows from the corresponding result [Io] for curves in the plane  since $M_{\theta}$ decouples into maximal averages over planes parallel to $P_{\theta}$. Since the estimates of [Io] are uniform
under small perturbations, a fact that derives from the corresponding uniformity under small perturbations of the estimates in the circular maximal
theorem, one sees that $M$ itself is bounded from $L^p$ to $L^p$ for $p > n$. Thus when combined with
Theorem 1.1a) and its smooth analogue described in section 5, we see that in the case where $\phi(x,y)$ is smooth and  the
tangent plane to $S$ at $(x_0,y_0,z_0)$ does not contain the origin,  without any restrictions on the Hessian determinant of $s(x,y)$ 
we have that $M$ is bounded on $L^p$ if $p > max({1 \over g},2)$. This is the main result of [IKeM].

\noindent {\bf Global extensions.}

If one wants to obtain a global theorem regarding maximal averages  of the form $(1.1)$ when $\psi$ is no longer assumed to be localized to
a neighborhood of a single point and when the singularities of 
$\psi$ are of the form considered in this paper, then one can use a partition of unity to
reduce the question to Theorem 1.1. Denote the support of $\psi$ by $A$. If $g_X$ denotes the index corresponding to $g$ at a point $X \in S$, then one obtains that if the Hessian 
determinant of the associated $s(x,y)$ are not identically zero, then  $M$ is bounded on $L^p$ for $p > \sup_{X \in A} \max({1 \over g_X}, 2) = 
\max({1 \over \inf_{X \in A} g_X}, 2)$.  The lower semicontinuity of $g_X$ implies that  that $ \inf_{X \in A} g_X$ is actually equal to $g_X$ for some
 $X \in A$. The examples used in the proof of part b) of Theorem 1.1 for such an $X$ then give a corresponding sharpness statement.

\noindent {\bf Extensions to smooth surfaces.}

Theorem 1.1 extends to the case of smooth surfaces that are not flat to infinite order at $(x_0,y_0,z_0)$, when $g$ is appropriately defined.
The analogue of the Hessian
 condition for the general smooth case is that the Hessian determinant of $s(x,y)$ not vanish to infinite order at $(x_0,y_0)$. Because the proof of this
extension involves a technical modification of the proof that might obscure the essence of the argument, we will not prove it in full detail. Instead, in section 5 we 
will provide a detailed sketch of the arguments. 

\subsection {Sobolev estimates for Radon transforms.}

We now come to our theorem concerning Sobolev space estimates for the Radon transform $R$. Let $T$ denote the interior of the region in the $xy$ plane bounded by the lines $y = x,\, y = 1 - x,\,y = 0$, and $y = g$. So  if $g \geq {1 \over 2}$ then $T$ is a triangle
whose upper vertex is $(1/2,1/2)$, and if $g < {1 \over 2}$ then $T$ is a trapezoid whose upper side is the portion of the line $y = g$ for which $x$ is
 in the interval $(g, 1- g)$. Our theorem is as follows.

\begin {theorem} There is a neighborhood $N$ of $(x_0,y_0)$ such that if 
$\phi(x,y)$ is supported in $N$ and $(1.3)$-$(1.4)$ are satisfied then we have the following.

\noindent {\bf a)} If  $({1 \over p},\alpha) \in T$, then $R$ is bounded from $L^p$ to $L^p_{\alpha}$.

\noindent {\bf b)} Suppose $g < 1$. If on some neighborhood of
 $(x_0,y_0,z_0)$ one has $|\phi(x,y)| > C|{\bf x} - {\bf x_0}|^{-\beta}$ for some $C > 0$, then $R$ is not bounded
 from $L^p$ to $L^p_{\alpha}$ for any $1 < p < \infty$ and $\alpha > g$.
 
\end{theorem}

So when $g < {1 \over 2}$ and ${1 \over p} \in (g, 1 - g)$, Theorem 1.2 
shows that the optimal $L^p$ Sobolev improvement (up to endpoints) that one can obtain is $g$ derivatives. A natural question to ask is if one can gain
$g$ derivatives for such $p$. It turns out that when $d$ in $(1.6)$ is equal to 1, then one does not gain $g$ derivatives, but when 
 $d = 0$ one gains $g$ derivatives for $p = 2$. The author does not know what happens when 
$p \neq 2$ in the $d = 0$ situation.

Analogous to the situation with the maximal averages, one can combine Theorem 1.2 with a straightforward partition of unity argument to 
prove an $L^p$ Sobolev improvement theorem for Radon transforms  when $\psi(X)$ in $(1.2)$ is not localized to near a specific point, when the density function has the types of singularities considered in this paper. 

Similarly to the situation of Theorem 1.1,
Theorem 1.2 extends to smooth surfaces via a simplified version of the arguments described in section 5 for the maximal averaging operators.

\subsection  {Some history.}

There has been a lot of work concerning local maximal averages over hypersurfaces with smooth density function $\phi(x,y)$. The initial work was
in [S1], where maximal averages over $n$-dimensional spheres were analyzed for $n \geq 2$ and $M$ was shown to be bounded on $L^p$ exactly 
when $p > {n + 1 \over n}$. The tricky case when $n = 1$ was later dealt with in [B], where boundedness of $M$ was shown to indeed hold if and only if
$p > 2$. These results can be generalized to situations where $S$ is a smooth hypersurface for which the Hessian determinant has positive rank, as was shown
in [So] and [Gr]. 
As for more general hypersurfaces, the paper [SoS] showed that if $S$ is a smooth hypersurface for which  the Gaussian curvature of $S$
does not vanish to infinite order near $(x_0,y_0,z_0)$, there is some $p < \infty$ for which $M$ is bounded on $L^p$. Optimal values of $p$ for which
$M$ is bounded on $L^p$ have been proven under a
nondegeneracy condition
on the Newton polyhedron [G4], as well as for convex hypersurfaces of finite line type such as in [CoMa] [IoSa] [NSeW].

For the $n = 2$ case that is being considered in this paper, the paper [IKeM] shows optimal $L^p$ boundedness of 
$M$ for all smooth $S$ when $\phi(x,y)$ is smooth, so long as the origin is not contained in the tangent plane to $S$ at $(x_0,y_0,z_0)$.
 Thus in the case of smooth $\phi(x,y)$, Theorem 1.1a) and its smooth extension in section 5 extend the sharp estimate of [IKeM] 
to the situation where the the tangent plane to $S$ at $(x_0,y_0,z_0)$ contains the origin unless the Taylor series of $s(x,y)$ 
vanishes to infinite order at the origin. The sharpness of this estimate when this tangent plane contains the origin is not understood in general, but as 
mentioned earlier, there are situations (see [Z] for more details) where one obtains stronger estimates. Also, as mentioned above,
the author earlier wrote a paper [G7] dealing with the $n=2$, smooth $\phi(x,y)$ situation.

As for Radon transforms and related operators, there has been a vast amount of work concerning boundedness properties between function spaces,
so we mainly restrict our attention to $L^p$ to $L^p_{\alpha}$ improvement results for hypersurfaces. For the case of translation-invariant smooth 
curves in $\R^2$, the sharp 
analogue to Theorem 1.2 is given by [Gra] and [C2]. For the more general non-translation-invariant situation for curves in $\R^2$, thorough $L^p_{\alpha}$ to
$L^q_{\beta}$ estimates up to endpoints are shown in [Se], which extend the $p = 2$ semi-translation-invariant situations that follow from [PS].

For the two dimensional translation-invariant situation of this paper, if $p = 2$ then the amount of Sobolev space
improvement is equal the exponent of decay of the associated surface measure Fourier transform. Hence when $\phi(x,y)$ is smooth and $p = 2$ the sharp estimates follow
from the analogous surface measure Fourier transform estimates of [IKeM] and [IM]. When $\phi(x,y)$ is singular, one correspondingly gets some sharp estimates for
the $p = 2$ case from the author's earlier work such as [G1]. It is worth pointing out that in many of the earlier results for smooth $\phi(x,y)$ such as
[D] and [IKeM], the proofs require $\phi(x,y)$ to be $C^3$, while the arguments of this paper only require that $\phi(x,y)$ be $C^1$ away from the 
origin.

When the exponent $\beta$ is close enough to 2, there is necessarily an interval $I$ containing 2 such that sharp
$L^p$ Sobolev improvement results for $p \in I$ follow from  [St]. In these cases it can be shown that the index $g$ of Theorems 1.1 and 1.2 
is necessarily equal to ${2 - \beta \over o}$, where $o$ is the order of the zero of $s(x,y)$ at $(0,0)$. If $\beta = 2$ and one assumes an appropriate cancellation condition on $\phi(x,y)$
near $(x_0,y_0)$, then the Radon transform becomes a singular Radon transform, and the general result of [CNSW] shows
$R$ is always bounded from each $L^p$ to itself for $1 < p < \infty$, as long as $s$ does not vanish to infinite order at $(x_0,y_0)$.

\subsection {Alternative formulations of the index $g$.}

There are a few alternative formulations of the index $g$ that will indicate the connection between $g$ and the corresponding indices in the statements of other
theorems on $L^p$ boundedness of maximal averages, especially those of [IKeM] and [IoSa]. First consider the case when $\phi(x,y)$ is smooth. By 
resolution of singularities (see [AGV] for details), if $(1.3)$ is satisfied and $s(x,y)$ is not identically zero,  there exists a positive number $o$, 
often called the oscillatory index of $s(x,y)$ at the origin, an integer 
$e = 0$ or $1$, and a neighborhood $N$ of the origin such that if $\phi(x,y)$ is supported in $N$ then as $\lambda \rightarrow \infty$ 
we have asymptotics of the form
$$ \int_{\R^2} e^{-i\lambda s(x,y)} \phi(x,y) \,dx\,dy = a_{\phi} \lambda^{-o}|\ln \lambda|^e + o(\lambda^{-o}|\ln \lambda|^e) \eqno (1.8)$$
Here $a_{\phi}$ is nonzero if $\phi(0,0) \neq 0$. Standard methods in resolution of singularities (again we refer to [AGV] for details) show that $o = g$, and also
in the notation of $(1.6)$ that  $e = d$ except possibly when the Hessian determinant of $s$ is nonvanishing at $(0,0)$. It turns out that the same
methods show that one has analogous asymptotics for the more general $\phi(x,y)$ satisfying $(1.4)$ considered in this paper, and that the generalization of
the statement that $o = g$ once again holds. Thus one could state Theorems 1.1 and 1.2 in terms of the oscillatory index $o$ if one prefers. 

One can connect the oscillatory index to Fourier transform decay estimates for surface measures by
looking for the optimal $q > 0$ for which there is an $f = 0$ or $1$ and a neighborhood $N$ of the origin such that for all $\phi$ supported on $N$
 one has an estimate of the form
$$\bigg|\int_{\R^2} e^{-i\lambda_1 s(x,y) - i\lambda_2 x - i\lambda_3 y} \phi(x,y)\,dx\,dy\bigg|  \leq C_{\phi}|\lambda|^{-q}|\ln \lambda|^f \eqno (1.9)$$
Here $|\lambda|$ denotes the magnitude of $|(\lambda_1,\lambda_2,\lambda_3)|$. Thus one is looking at how the upper bounds in $(1.8)$ change
if one requires that they still hold after a linear perturbation of the phase. Equivalently, $q$ is the supremum of the $q'$ for which the Radon transform $R$ gains
$q'$ derivatives on $L^2$. Thus by Theorem 1.2, when $g \leq {1 \over 2}$ one has $q = g$, and when $g > {1 \over 2}$ one has
$q \geq {1 \over 2}$. Hence in the  $g \leq {1 \over 2}$ situation  Theorems  1.1 and 1.2 could be stated in terms of $q$.

One could go further and not only take the infimum of the oscillatory index under linear perturbations of the phase but also over $(x,y)$ in a neighborhood
of $(x_0,y_0)$, as was done in [IKeM]. Because $(1.6)$ holds for all $\phi$ supported in a neighborhood $N$ of the origin, the index analogous to $g$
for all points in $N$ must be at least $g$. In other words, the index $g$ is a lower semicontinuous function. Similar considerations show the oscillatory index is
lower semicontinuous. Hence whenever $g \leq {1 \over 2}$ the infimum
of the oscillatory index over all all linear perturbations of the phase and all $(x,y) \in N$ must once again be $g$.

Next, rewrite $(1.6)$ as follows, where $X$ denotes $(x,y,z)$, $X_0$ denotes $(x_0,y_0,z_0)$, and ${\bf n}$ denotes the normal direction to $S$ at $(x_0,y_0,z_0)$.
$$\int_{\{X \in S:\,|{\bf x} - {\bf x_0}| < r,\,|(X - X_0) \cdot {\bf n}| < \epsilon\}} |{\bf x} - {\bf x_0}|^{-\beta} \,d {\bf x} = C_{\beta, r} \epsilon^{g}
|\ln \epsilon|^{d} + o(\epsilon^{g}
|\ln \epsilon|^{d}) \eqno (1.10)$$
One can replace ${\bf n}$ by another direction and ask how the asymptotics $(1.10)$ change. It is not hard to show that the index $g$ increases.
For in any other direction, the integral $(1.10)$ will effectively be an integral over a slat of width $C\epsilon$ centered at $(0,0)$ and by performing 
the calculation one obtains an upper bound of $C\epsilon^{\min(2 - \beta, 1)}$. On the other hand, the growth rate in $(1.10)$ will be at fastest
the growth rate for the case when $s(x,y) = x^2 + y^2$, and again doing a calculation one sees in this case that $g = 1 - {\beta \over 2} <
\min(2 - \beta, 1)$. Thus the normal direction always gives the slowest growth rate in $\epsilon$. Hence one could have defined $g$ as the infimum of 
the growth rate exponents over all directions. 

One can pursue this idea further and not only take the infimum of the growth rate exponents over all directions but also over all $(x,y,z)$ in a neighborhood
of $(x_0,y_0,z_0)$ (not even restricting to points on $S$). Call this infimum $h$. By the lower semicontinuity of $g$, 
one has that $h = g$. Thus the statements of Theorems 1.1 and 1.2 could have been
 reformulated in terms of $h$ in place of $g$. This type of formulation is made in [IoSa], where they conjecture that if $S$ is a smooth 
compact surface in any dimension, and $\phi(x,y) =1 $,  then if 
$p > 2$ the operator $M$ is bounded on $L^p$ if and only if $|d(S,H)|^{-{1 \over p}}$ locally has finite integral for any hyperplane $H$ not passing through the origin, 
where $d(S,H)$ denotes the distance from $S$ to the hyperplane $H$. 

\section {The resolution of singularities theorems and some consequences.}

We will make use of a couple of resolution of singularities theorems from the author's earlier work. For the first, we first rotate coordinates
so that $\partial_x^o s(0,0) \neq 0$ and $\partial_y^o s(0,0) \neq 0$, where $o \geq 2$ is order of the zero of $s(x,y)$ at $(0,0)$. Let
$H(x,y)$ denote the Hessian determinant of $s(x,y)$. If $H(0,0) = 0$ we also 
assume the rotation is such that $\partial_x^p H(0,0) \neq 0$ and $\partial_y^p H(0,0) \neq 0$ for some $p \geq 1$. We then apply 
Theorem 2.1 of [G2] to the rotated $s(x,y)$, which gives the following.

Write the Taylor expansion of $s(x,y)$ at the origin as $\sum_{\alpha, \beta} s_{\alpha, \beta} x^{\alpha} y^{\beta}$. Divide the
$xy$ plane into eight triangles by slicing the plane using  the $x$ and $y$ axes  and two lines through the origin, one of the form $y = mx$ for some $m > 0$ and one of the form $y = mx$ for some $m < 0$. One must ensure that these two lines are not ones  on 
which the function $\sum_{\alpha + \beta = o} s_{\alpha, \beta}x^{\alpha}y^{\beta}$ vanishes other than at the
origin.  After reflecting about the $x$ and/or $y$ axes and/or the line $y = x$ if necessary, each of the triangles becomes of the form $T_b = \{(x,y) \in \R^2: x > 0,\,0 < y < bx\}$ (modulo an inconsequential boundary set of measure zero). The version of the real analytic case of Theorem 2.1 of [G2] that is pertinent 
here is what was called Theorem 2.1 in [G1]:

\begin {theorem} (Theorem 2.1 of [G1]) Let  $T_b = \{(x,y) \in \R^2: x > 0,\,0 < y < bx\}$ be as above. Abusing notation slightly, use the notation $s(x,y)$ to denote the reflected function $s(\pm x,\pm y)$ or $s(\pm y, \pm x)$ corresponding to $T_b$.
 Then there is a $a > 0$ and a positive integer $N$ such that
if $F_a$ denotes  $\{(x,y) \in \R^2: 0 \leq x\leq a, \,0 \leq y \leq bx\}$, then one can write $F_a = \cup_{i=1}^n cl(D_i)$, such that for to each $i$ there is a $k_i(x) = p_i x$ or $k_i(x) = p_i x + l_i x^{s_i} + ...$ with $k_i(x^N)$ real analytic, $l_i \neq 0$, and $s_i > 1$, such that after a coordinate change of the form $\eta_i(x,y) = (x, \pm y + k_i(x))$, the set $D_i$ becomes a set $D_i'$ on which the function $s \circ \eta_i(x,y)$ approximately becomes a monomial $d_i x^{\alpha_i}y^{\beta_i}$, $\alpha_i$ a nonnegative rational number and $\beta_i$ a nonnegative integer in the following sense.

\noindent {\bf a)} $D_i' = \{(x,y): 0 < x < a, \, g_i(x) < y < G_i(x)\}$, where $g_i(x^N)$ and $G_i(x^N)$ are
real analytic. If we expand $G_i(x) =  H_i x^{M_i} + ...$, then $M_i \geq 1$ and $H_i > 0$.

\noindent {\bf b)} Suppose $\beta_i = 0$. Then $g_i(x) = 0$. The set $D_i'$ can
be constructed such that for any preselected $\eta > 0$ there is a $d_i \neq 0$ such that on $D_i'$, for all $0 \leq l \leq \alpha_i$ one has
$$ |\partial_x^l (s \circ \eta_i)(x,y) -  d_i\alpha_i(\alpha_i -1) ... (\alpha_i - l + 1)x^{\alpha_i - l}| < \eta|d_i| x^{\alpha_i-l} \eqno (2.1)$$ 

\noindent {\bf c)} If $\beta_i > 0$, then $g_i(x)$ is either identically zero or $g_i(x)$ 
can be expanded as $h_ix^{m_i} + ...$ where $h_i > 0$ and $m_i > M_i$. The $D_i'$ can
be constructed such that such that for any preselected $\eta > 0$ there is a $d_i \neq 0$ such that on $D_i'$, for all $0 \leq l \leq \alpha_i$ and all $0 \leq  m \leq \beta_i$ one has
$$|\partial_x^l\partial_y^m(s \circ \eta_i)(x,y) -  \alpha_i(\alpha_i - 1)...(\alpha_i - l + 1)\beta_i(\beta_i - 1)...(\beta_i - m + 1)
d_ix^{\alpha_i - l}y^{\beta_i - m}| $$
$$\leq \eta |d_i| x^{\alpha_i - l}y^{\beta_i - m} \eqno (2.2)$$
\end {theorem}

In [G2] it was shown that one can may do the constructions so that $s_i \leq M_i$ for all $i$ whenever $k_i(x)$ is not of the form
$p_i x$.

On the domains $D_i'$ for which $\beta_i = 0$ and $k_i(x)$ is not of the form $p_i x$ (i.e. there exists a nonzero $l_i x^{s_i}$ term) we must do
a  second resolution
of singularities, this time simultaneously resolving the singularities of $s \circ \eta_i(x,y)$, $\partial_y (s \circ \eta_i(x,y))$, and $\partial_{yy} (s \circ \eta_i(x,y))$. (Because the coordinate changes are effectively translations in $y$ for fixed $x$, such $y$ derivatives commute with the coordinate changes and 
thus such a resolution of singularities makes sense.) To perform this simultaneous resolution of singularities, we use following theorem from [G3]. Although
it is stated for real analytic functions of $x$ and $y$, the same proof holds for real analytic functions of $x^{1 \over N}$ and $y$ for a positive integer $N$.

\begin {theorem} (Theorem 2.2 of [G3]) Suppose $S_1(x,y),...,S_k(x,y)$ are real analytic functions on a neighborhood of the origin with
$S_j(0,0) = 0$ for each $j$. Let $D_i'$, $\alpha_i$, and $\beta_i$ be as in Theorem 2.1 applied to  $\prod_{j=1}^k S_j(x,y)$.
 Then one can further divide each $D_i'$ into finitely many pieces $D_{il}$, such that on each $D_{il}$ an
additional coordinate change of the form $(x,y) \rightarrow (x,y - c_{il} x^{M_i})$ or $(x,y -  c_{il}x^{m_i})$, $c_{il} \geq 0$, will
result in each $S_j(x,y)$ satisfying the conclusions of Theorem 2.1, with one difference: Let the domains in the new coordinates be denoted by
$D_{il}'$. Then the $D_{il}'$ 
can now only be assumed to have the same form as the domains where $\beta_i > 0$ in Theorem 2.1. That is, $D_{il}'$ has the form $\{(x,y): 0 < x < a, \, g_{il}(x) < y < G_{il}(x)\}$, where $g_{il}(x^K)$ and $G_{il}(x^K)$ are real analytic for some positive integer $K$,  $G_{il}(x) =  H_{il} x^{M_{il}} + ...$, and  $g_{il}(x)$ is identically zero or is of the form $h_{il} x^{m_{il}} + ... $ where $1 \leq M_{il} < m_{il}$ and $h_{il}, H_{il} > 0$.
\end{theorem}

\noindent We will also make use of the following corollary to Theorem 2.2, which follows from Corollary 2.3 of [G3].

\begin {corollary} Let $\eta_{il}(x,y)$ denote the composition of the coordinate changes in Theorem 2.2. For any given $K$, however large, the $D_{il}'$ can be constructed 
so that there is a constant $C_K$ so that on $D_{il}'$ one has $|\partial_x^{a}\partial_y^{b} (S_j \circ\eta_{il}) (x,y)| \leq C_K x^{-a}y^{-b} |S_j \circ \eta_{il}(x,y)|$ 
for all $a, b < K$.
\end{corollary}

We now describe some consequences of Theorems 2.1 and 2.2 that we will make use of in section 4. These facts are best described in terms of Newton 
polygons. Let $f(x,y)$ denote a power series in $x^{1 \over N}$ and $y$ for some positive integer $N$, and write $f(x,y) = \sum_{a,b} f_{a,b}
x^a y^b$.

\begin {definition} For any $(a,b)$ for which $f_{ab} \neq 0$, let $Q_{ab}$ be the
quadrant $\{(x,y) \in \R^2: 
x \geq a, y \geq b \}$. Then the {\it Newton polygon} $N(f)$ of $f(x,y)$ is defined to be 
the convex hull of the union of all $Q_{ab}$.  
\end {definition}

 The boundary of $N(f)$ consists of finitely many (possibly none) bounded edges of negative slope
as well as an unbounded vertical ray and an unbounded horizontal ray. We write these slopes as $-{1 \over p_n} < -{1 \over p_{n-1}} < ...
< -{1 \over p_0}$, where $p_n = 0$ and $p_0 = \infty$.  We denote by $(a_i,b_i)$ the vertex of $N(f)$ joining the edge of slope $-{1 \over p_i}$ 
to the edge of slope $-{1 \over p_{i+1}}$. 

We focus on the case where $f(x,y)$ is a real analytic function of 
$x^{1 \over N}$ and $y$ on a neighborhood of the origin. If $0 < i < n-1$ there exists a
constant  $c_i$ such that for each $M > 0$ there is a $\delta > 0$ such that if ${1 \over \delta}|x|^{p_i} < |y| < \delta |x|^{p_{i+1}}$ then 
${1 \over 2} < |{f(x,y)  \over c_i x^{a_i}y^{b_i}}| < 2$ and $| {x^{a_j}y^{b_j} \over  x^{a_i}y^{b_i}}| \leq {1 \over M}$ for $j \neq i$. When $i = 0$
 the same is true if we replace the condition ${1 \over \delta}|x|^{p_i} < |y| < \delta |x|^{p_{i+1}}$ by the condition that $|y| < \delta |x|^{p_1}$,
and if $i = n-1$ the same is true  if we replace the condition ${1 \over \delta}|x|^{p_i} < |y| < \delta |x|^{p_{i+1}}$ by the condition that $|y| 
> {1 \over \delta} |x|^{p_{n-1}}$. For brevity, we refer to the proof of Lemma 2.4 of [G5] for a proof of of a slight variant of these facts (the analogue for
smooth functions) rather than present the full argument here.

We now give some pertinent consequences of the above considerations. If $(\alpha_i,\beta_i)$ is as in Theorem 2.1 or 2.2 for $s(x,y)$ or some
 $S_j(x,y)$ respectively, then $(\alpha_i,\beta_i)$ is a vertex of the Newton polygon of the function. Furthermore, if $(\alpha_i,\beta_i)$ is between the edges
 with slopes $-{1 \over p_{k+1}}$ and $-{1 \over p_k}$, then 
the numbers called $M_i$ or $M_{il}$ in Theorem 2.1 and 2.2 respectively must satisfy $M_i \geq p_{k + 1}$ or $M_{il} \geq p_{k + 1}$. If
$h_i$ or $h_{il}$ is nonzero, then one similarly has $m_i \leq p_k$ or $m_{il}\leq p_k$. 

Recall we are first applying Theorem 2.1 to $s(x,y)$ and then we are applying
Theorem 2.2 to  $s \circ \eta_i(x,y)$, $\partial_y (s \circ \eta_i(x,y))$, and $\partial_{yy} (s \circ \eta_i(x,y))$ on $D_i'$ in the cases
where $\beta_i = 0$ and $k_i(x)$ is not linear after the application of Theorem 2.1. Thus if Theorem 2.2 is applied,  after the application of Theorem 2.1 
$s \circ \eta_i(x,y)$ was comparable to $x^{\alpha_i}$ and the lower edge of $D_i'$ was on the $x$-axis. As a result, $(\alpha_i, 0)$ was the lowest vertex of the Newton polygon of $s \circ \eta_i(x,y)$ after this application of Theorem 2.1. Consequently, 
if $\zeta_{il}(x,y)$ denotes the coordinate change analogous to $\eta_i(x,y)$
for this application of Theorem 2.2, then  $s_{il}(x,y) = s \circ \eta_i \circ \zeta_{il} (x,y)$ is still comparable to $x^{\alpha_i}$ after applying Theorem 2.2 and $(\alpha_i,0)$ is still the lowest vertex of the Newton polygon of $s_{il}(x,y)$.

Next, let $(\gamma_{il},\delta_{il})$ be such that
$\partial_{yy} s_{il}(x,y)$ is comparable to $x^{\gamma_{il}}y^{\delta_{il}}$ after the above application of Theorem 2.2. 
Let $(A_{il},B_{il})$ denote the upper vertex of the lowest non-horizontal edge of the Newton polygon of $s_{il}(x,y)$
and let $(0,B_{il}')$ denote the intersection of the $y$ axis with the line extending this edge. 

\begin {lemma} If $B_{il} \geq 2$, then on the domain $D_{il}'$ there are constants $C$ and $C'$ such that  then 
$$y^{B_{il}'} \leq C x^{A_{il}}y^{B_{il}} \leq C'x^{\gamma_{il}}y^{\delta_{il} + 2} \eqno (2.3)$$
\end {lemma}

\noindent {\bf Proof.} By the above discussion, $(\gamma_{ij},\delta_{ij})$ is a vertex of the Newton polygon of $\partial_{yy} s_{il}(x,y)$, and if $(v_1,v_2)$ is any other vertex, on the domain $D_{il}'$ one has that $x^{v_1}y^{v_2} < 
Cx^{\gamma_{il}}y^{\delta_{il}}$ for some constant $C$. Since $B_{il} \geq 2$, one has that $(A_{il},B_{il} - 2)$ is a vertex of the Newton polygon of $\partial_{yy} s_{il}(x,y)$.
Hence one may take $v_i = A_{il}$ and $v_2 = B_{il} - 2$ and the right-hand inequality of $(2.3)$ follows. 

As for the left-hand inequality, the lowest non-horizontal edge of $N(f)$ for $f = s_{il}(x,y)$ connects $(A_{il},B_{il})$ to 
$(\alpha_i,0)$, and since the lowest vertex dominates here we have $x^{A_{il}}y^{B_{il}} < Cx^{\alpha_i}$ on $D_{il}'$. Since $(0,B_{il}')$ is on the line 
containing this edge,  for some $e > 0$ we have 
$$(y^{B_{il}'} / x^{A_{il}}y^{B_{il}}) = (x^{A_{il}}y^{B_{il}} / x^{\alpha_i})^e < C^e \eqno (2.4)$$
This gives the left-hand inequality of $(2.3)$ and we are done.

\begin{lemma}
 Let $s_{il}^*(x,y) = |y^{B_{il}'} s_{il}(x,y)|^{1 \over 2}$. Then $\int_{D_{il}'}(s_{il}^*(x,y))^{-t}|(x,y)|^{-{\beta}}\,dx\,dy$
 is finite for all $t < g$, where $g$ as in Theorems 1.1 and 1.2.
 \end{lemma}

\noindent {\bf Proof.} Let $n_{il}$ be such that the slope of the lowest nonhorizontal edge of the Newton polygon of  $s_{il}(x,y)$
has slope $-{1 \over n_{il}}$. So $|s_{il}(x,y)| \sim x^{\alpha_i}$ on a set of the form $0 < x < r$, $0 < y < \nu x^{n_{il}}$ for some small
$r$ and $\nu$. The definition of $g$ implies that $\int |s(x,y)|^{-t}|(x,y)|^{-\beta}$ is finite on a neighborhood of the origin if $t < g$. As a result, we 
have
$$\int_{\{(x,y): 0 < x < r,\, 0 < y < \nu x^{n_{il}}\}}|s_{il}(x,y)|^{-t}|(x,y)|^{-\beta}\,dx\,dy < \infty \eqno (2.5)$$
Since $|s_{il}(x,y)| \sim x^{\alpha_i}$ and $|(x,y)| \sim x$ on this domain of integration, 
whenever $t < g$ we have that 
$$ \int_{\{(x,y): 0 < x < r,\, 0 < y < \nu x^{n_{il}}\}}x^{-\alpha_i t - \beta}\,dx\,dy < \infty \eqno (2.6)$$
As a result, for all $t < g$ one has $t < {1 + n_{il} -\beta \over \alpha_i}$. In other words, $g \leq {1 + n_{il} - \beta \over \alpha_i}$. Next, note that 
the definition of $B_{il}'$ can be rewritten as just ${\alpha_i \over n_{il}}$. In addition, since $x^{\alpha_i}$ dominates the Taylor expansion of 
$s_{il}(x,y)$ on $D_{il}'$, one has that $D_{il}' \subset \{(x,y):\, 0 < x < 1,\,0 < y < Cx^{n_{il}}\}$ for some $C$. So we have
$$\int_{D_{il}'} (s_{il}^*(x,y))^{-t}|(x,y)|^{-\beta}\,dx\,dy = \int_{D_{il}'} |y^{B_{il}'} s_{il}(x,y)|^{-{t \over 2}}|(x,y)|^{-\beta}\,dx\,dy$$
$$\leq  C\int_{\{(x,y):\, 0 < x < 1,\,0 < y < x^{n_{il}}\}} y^{-{\alpha_i t \over 2n_{il}}}x^{-{\alpha_i t\over 2}}x^{-\beta}\,dx\,dy \eqno (2.7)$$
Performing the $y$ integration first in $(2.7)$, we see that $(2.7)$ is finite if ${\alpha_i t \over 2n_{il}} < 1$ and $t < {1 + n_{il} - \beta \over \alpha_i}$. The latter
condition holds due to the finiteness of $(2.6)$, and the former condition holds since $t < {1 + n_{il} - \beta \over \alpha_i} \leq {2n_{il} \over \alpha_i}$.
This completes the proof of Lemma 2.6.

\section {Preliminary lemmas and an overview of the proofs of Theorems 1.1 and 1.2.}

\subsection {Preliminary lemmas}

It is well-known in the field (we refer to chapter 11 of [S2] for details)  that complex interpolation between $L^2$ and $L^{\infty}$ bounds for damped
versions of a given maximal average can often be used in proving optimal $L^p$
boundedness of the original maximal average. As was described in [SoS], in such an interpolation the following lemma is useful. 
It provides a way of reducing $L^2$  boundedness of maximal averages to oscillatory integral decay estimates, and has been used in various papers in
this subject, including [IoSa].

\begin{theorem} ([SoS]) Suppose $d\sigma$ is a measure on $\R^n$ such that for some $\nu > 0$ the following holds for all 
multiindices with $|\alpha| = 0, 1$.
$$|\partial_{\lambda}^{\alpha}\widehat{d\sigma}(\lambda)| \leq A(1 + |\lambda|)^{-{1 \over 2} - \nu} \eqno (3.1)$$
Let $M_{\sigma}$ denote the maximal operator 
$$M_{\sigma} f(x) = \sup_{t > 0} \bigg|\int_{\R^n} f(x - ts) d\sigma(s)\bigg| \eqno (3.2)$$
Then there is a constant $B$ depending on $A$ and $\nu$ such that 
$$||M_{\sigma} f||_{L^2} \leq B||f||_{L^2}\eqno (3.3)$$
\end{theorem}

In order to prove the needed Fourier transform decay estimates in our setting, we will make use of two Van der Corput style theorems. The first is the 
standard Van der Corput lemma (see p. 334 of [S2]).

\begin{lemma} Suppose $P(x)$ is a real-valued $C^k$ function on the interval $[a,b]$ with $|P^{(k)}(x)| > M$ on $[a,b]$ for
some $M > 0$. Let $\psi(x)$ be a complex-valued $C^1$ function on $[a,b]$. If $k \geq 2$ there is a constant $c_k$ depending only on $k$ such that
$$\bigg|\int_a^b e^{iP(x)}\psi(x)\,dx\bigg| \leq c_kM^{-{1 \over k}}\bigg(|\psi(b)| + \int_a^b |\psi'(x)|\,dx\bigg) \eqno (3.4)$$
If $k =1$, the same is true if we add the conditions that $P(x)$ is $C^2$ and that $P'(x)$ is monotonic on $[a,b]$. 
\end{lemma}

The second Van der Corput style lemma we will use is a version proved in [G1] that holds for mixed partial derivatives.

\begin{lemma} Let $I_1$ and $I_2$ be closed intervals of lengths $l_1$ and $l_2$ respectively, and for some strictly
monotone functions $f_1(x)$ and $f_2(x)$ on $I_1$ with $f_1(x) \leq f_2(x)$ let $R = \{(x,y) \in I_1 \times I_2: f_1(x) \leq y \leq  f_2(x)\}$ (Note $R$ might just be $I_1 \times I_2$). Suppose for some $k \geq 2$, $P(x,y)$ is a $C^k$ real-valued function on $R$ such that for each $(x,y) \in R$ one has
$$|\partial_{xy} P(x,y)| > M\,\,\,\,\,\,\,\,\,\,\,\,\,\,\,\,{\rm and }\,\,\,\,\,\,\,\,\,\,\,\,\,\,\,
\partial_y^k P(x,y) \neq 0 \eqno (3.5a)$$
Further suppose that $\Psi(x,y)$ is a complex-valued function on $R$ that is $C^1$ in the $y$ variable for fixed $x$,  such that 
$$ |\Psi(x,y)| < N \,\,\,\,\forall x,y\,\,\,\,\,\,\,\,\,\,\,\,\,\,\,{\rm and }\,\,\,\,\,\,\,\,\,\,\,\,\,\,\,\,\int_{\{y: (x,y) \in R\}} |\partial_y\Psi(x,y)|\,dy< N \,\,\,\,\forall x \eqno (3.5b)$$
If  $R' \subset R$ such that the intersection of $R'$ with each vertical line is either empty or is a set of at most $l$ intervals, then 
the following estimate holds. 
$$\bigg|\int_{R'} e^{i P(x,y)}\Psi(x,y)\,dx\,dy\bigg| < C_{kl}  N \bigg({l_1l_2 \over M}\bigg)^{1 \over 2} \eqno (3.6)$$
\end{lemma}

The following lemma characterizes when the Hessian determinant of a smooth function near the origin has an identically zero Taylor series at the origin.

\begin{lemma} Let $f(x,y)$ be a smooth function on a neighborhood of the origin with a zero of order $n > 1$ at the origin, such that
the Hessian determinant of $f(x,y)$ vanishes to infinite order at the origin. Then 
there is a linear map $L$ such that the Taylor series of $f(L(x,y))$ at the origin is of the form $a(x,y)y^n$, where $a(0,0) \neq 0$.
\end{lemma}

\noindent {\bf Proof.} Let $f(x,y)$ be any smooth function on a neighborhood of the origin with a zero of order $n$ at the origin. Write the Taylor 
expansion of $f(x,y)$ at the origin in the form $\sum_{a,b} f_{a,b} x^a y^b$, and the Taylor expansion of
the Hessian determinant of $f(x,y)$ at the origin as  $\sum_{a,b} h_{a, b} x^a y^b$. For any $m > 0$, let $e_m$ be the minimum value 
of $a + mb$ amongst all nonzero $f_{a,b}$, and let $f_m(x,y) = \sum_{a + mb  = e_m} f_{a,b} x^a y^b$. Similarly, if the Taylor expansion of the Hessian
is nonzero let $l_m$ be  the minimum value 
of $a + mb$ amongst all nonzero $h_{a,b}$, and let $h_m(x,y) = \sum_{a + mb  = l_m} h_{a,b} x^a y^b$. If the Hessian determinant of $f_m(x,y)$ is
not identically zero, then the Taylor expansion of $h(x,y)$ at the origin is equal to the Hessian determinant of $f_m(x,y)$ plus possibly some terms for which $a + mb$ is
 geater than $l_m$. Thus in any situation where the Taylor expansion of $h(x,y)$ vanishes to infinite order at the origin, the Hessian determinant of each
 $f_m(x,y)$ must be also be identically zero.
 
 Now we suppose that the Taylor expansion of $h(x,y)$ vanishes to infinite order at the origin. If the Newton polygon of $f$ had a vertex $(\alpha,\beta)$
 with $\alpha, \beta \neq 0$, then we could find some $m$ for which $f_m(x,y)$ is a multiple of $x^{\alpha}y^{\beta}$, a function whose Hessian determinant is of the form
 $cx^{2\alpha -2}y^{2\beta - 2}$ for some $c \neq 0$  and is therefore not identically zero. Thus the only vertices of the Newton polygon lie on the $x$ and $y$ axes. If there
 is a vertex on only one of these two axes, then we are done. Otherwise, if
 $-{1 \over m}$ denotes the slope of the edge of the Newton polygon connecting the two vertices, then $f_m(x,y)$ must have vanishing Hessian determinant.
  But
 some algebra shows (see Corollary 2.2 of [dBvdE] for a more general result) that $f_m(x,y)$ must be of the form $(cx + dy)^n$ for some nonzero 
 $c$ and $d$. So if
 we do the linear coordinate $L$ change turning $(x,y)$ to $(x,y - (c/d) x)$, then $f(L(x,y))$'s Newton polygon has a vertex at $(0,n)$ but not at
 $(n,0)$. Since the Hessian of $f(L(x,y))$ is also identically zero, by the above considerations this new Newton polygon can have vertices only on the
 $x$ or $y$ axes. If it had vertices on both axes, then like before the vertices would have to be $(0,n)$ and $(n,0)$. Since there is no vertex at $(n,0)$, we
 conclude that the Newton polygon has exactly one vertex, at $(n,0)$. This completes the proof of Lemma 3.4.

\subsection  {Overview of the proof of Theorem 1.1.}

We rewrite $(1.1)$ in coordinates for which $(0,0,1)$ is normal to $S$ at $(x_0,y_0,z_0)$. Let ${\bf x_0} = (x_0,y_0)$ and 
${\bf x'} = (x',y')$. We then have
$$M f(X) = \sup_{t > 0} \bigg| \int f(x - tx', y - ty', z - t (z_0 + s({\bf x'} - {\bf x_0})))\, \phi(x',y')\, dx ' dy'  \bigg|  \eqno (3.7)$$

For our complex interpolation, we embed the maximal operator $M$ in an analytic family of operators as follows. Let $H(x,y)$ denote the Hessian determinant of 
$s(x,y)$. Then for a small value of $\delta$ to be determined by our arguments, we look at the operators $M_z$ defined as follows. Recall that
$(A_{il},B_{il})$ denotes the upper vertex of the lowest edge of the Newton polygon of $s_{il}(x,y)$ after applying Theorem 2.2. Let $\bar{s}(x,y)$ denote
the function which is equal to $s(x,y)$ everywhere except on those domains of section 2 where $\beta_i = 0$, where $\eta_i(x,y)$ was nonlinear so that we had
to do a second resolution of singularities, and where $B_{il} \geq 2$. On the domains where these three conditions hold, one defines $\bar{s}(x,y)$ in the coordinates of $D_{il}'$
 to be the function $s_{il}^*(x,y) = |y^{B_{il}'}s_{il}(x,y)|^{1 \over 2}$ of Lemma 2.6. We then define $M_z$ by
$$M_z f(X) =  \sup_{t > 0} \bigg| \int f(x - tx', y - ty', z - t (z_0 + s({\bf x'} - {\bf x_0})) \, |\bar{s}({\bf x'} - {\bf x_0})|^z $$
$$\times |H({\bf x'} - {\bf x_0})|^{\delta z} e^{z^2} \phi(x',y')\, dx ' dy'  \bigg|  \eqno (3.8)$$
We will see in the notation of Theorem 1.1 that if $s > -g$ is fixed, then  if $\delta$ is sufficiently small one has the estimate 
$||M_z f||_{\infty} \leq C ||f||_{\infty}$, where $C$ is uniform over all $z$ with $Re\, z = s$. This will follow relatively easily by simply showing that
the measures of the surfaces in $(3.8)$ are uniformly bounded over such $z$. The vast majority of our effort will go into showing that if 
$s > max (0, {1 \over 2} - g)$ is fixed, then as long as the Hessian determinant of $s(x,y)$ is not identically
 zero, if $\delta$ is sufficiently small one has estimates $||M_z f||_2 \leq C ||f||_2$ for a constant $C$ that is uniform over all $z$ with $Re\, z = s$.
Theorem 3.1 will reduce these $L^2$ estimates to proving a Fourier transform decay estimate which will be the bulk of our effort.
Using complex interpolation will then give us Theorem 1.1.

To give an idea of how the Fourier transform decay estimates are proved, let $d\sigma_z$ denote the surface measure being dilated in $(3.8)$. 
Then shifting coordinates to be centered at $X_0 = (x_0,y_0,z_0)$ we have
$$\widehat{d\sigma_z}(\lambda) = e^{-i\lambda \cdot X_0}\int_{\R^2} e^{-i\lambda_1 x - i\lambda_2y - i\lambda_3s(x,y)} 
|\bar{s}(x,y)|^z |H(x,y)|^{\delta z}e^{z^2} \phi(x + x_0,y + y_0)\,dx\,dy \eqno (3.9)$$
Note that any $\lambda_i$ derivative of $(3.9)$ is of the same form as $(3.9)$ except with $\phi(x + x_0,y + y_0)$ replaced by $\phi(x + x_0,y + y_0)$ times a
smooth function. Thus our arguments bounding $|\widehat{d\sigma_z}(\lambda)|$ will always lead to the same estimates for each
$|\partial_{\lambda_i}\widehat{d\sigma_z}(\lambda)|$. Thus for the purposes of applying Theorem 3.1 in this paper, we will always focus on 
bounding $|\widehat{d\sigma_z}(\lambda)|$, with the understanding that the same argument will always give the same bound for each $|\partial_{\lambda_i}\widehat{d\sigma_z}(\lambda)|$.

Next, note that it makes sense to define $\rho(x,y) = \phi(x + x_0,y + y_0)$ and rewrite $(3.9)$ as 
$$\widehat{d\sigma_z}(\lambda) = e^{-i\lambda \cdot X_0}\int_{\R^2} e^{-i\lambda_1 x - i\lambda_2y - i\lambda_3s(x,y)} 
|\bar{s}(x,y)|^z |H(x,y)|^{\delta z}e^{z^2} \rho(x,y)\,dx\,dy \eqno (3.9')$$
The conditions $(1.4)$ become
$$|\rho(x,y)| \leq A|(x,y)|^{-\beta}{\hskip 0.7 in} |\nabla \rho(x,y)| \leq A|(x,y)|^{-\beta - 1} \eqno (3.10)$$
Let $P_{\lambda}(x,y) = \lambda_1 x + \lambda_2y + \lambda_3s(x,y)$ denote the phase function in $(3.9')$. When $|\lambda_1| + |\lambda_2| > |\lambda_3|$, 
one has that $|\nabla P_{\lambda}(x,y)| > {1 \over 4}|\lambda|$ on a sufficiently small neighborhood, and we will see that applying Van der Corput's lemma for
first derivatives on $(3.9')$ appropriately will give the estimates needed to apply Theorem 3.1. 

Thus the main effort is the situation where $|\lambda_1| + |\lambda_2| \leq |\lambda_3|$. For this, the first step will be to divide $(3.9')$ into two
pieces, depending whether or not $|H(x,y)| > |\lambda|^{-{1 \over 100}}$ and $|(x,y)| > |\lambda|^{-{1 \over 100}}$. 
Specifically we write $\widehat{d\sigma_z}(\lambda) = I_1(\lambda) + 
I_2(\lambda)$, where
$$I_1(\lambda) =  e^{-i\lambda \cdot X_0}\int_{\{(x,y):\, |H(x,y)| > |\lambda|^{-{1 \over 100}},\, |(x,y)| > |\lambda|^{-{1 \over 100}}\}} e^{-i\lambda_1 x - i\lambda_2y - i\lambda_3s(x,y)} |\bar{s}(x,y)|^z$$
$$\times  |H(x,y)|^{\delta z} e^{z^2} \rho(x,y)\,dx\,dy \eqno (3.11a)$$
$$I_2(\lambda) =  e^{-i\lambda \cdot X_0}\int_{\{(x,y):\, |H(x,y)| < |\lambda|^{-{1 \over 100}}{\,\, {\rm or}\,\,}|(x,y)| < |\lambda|^{-{1 \over 100}}\}} e^{-i\lambda_1 x - i\lambda_2y - i\lambda_3s(x,y)} |\bar{s}(x,y)|^z $$
$$ \times |H(x,y)|^{\delta z} e^{z^2} \rho(x,y)\,dx\,dy \eqno (3.11b)$$
On the domain of the first integrand, the Hessian determinant of $P_{\lambda}(x,y)$, given by $\lambda_3 H(x,y)$, will be of absolute value at
least $C|\lambda|^{99 \over 100}$.
Because this determinant is so large, we will see that using an elaboration of a type of argument often used for nondegenerate phases will provide the 
needed estimates to apply Theorem 3.1.

For $I_2(\lambda)$, we will need to delve deeper. After applying the resolution singularities algorithm of the last section to $s(x,y)$ and its first two $y$ 
derivatives as described above $(2.5)$, we have
domains $D_i$ and $D_i'$ if only Theorem 2.1 is being used, and domains $D_{il}$ and $D_{il}'$ if Theorem 2.1 followed by Theorem 2.2 on $D_i'$ are
being used. Recall the latter occurs when after applying Theorem 2.1 to $s(x,y)$ one has $\beta_i = 0$ and the coordinate change 
$\eta_i$ is not linear. To make
notation consistent, in the situation where we are applying Theorem 2.1 and Theorem 2.2 we write the combined coordinate change which we 
called $\eta_i \circ \zeta_{il}$ before as simply $\eta_{il}$. 

We will separately bound the contribution to $(3.11b)$ arising from each domain $D_i$ and $D_{il}$. We perform the coordinate change $\eta_i$ or
$\eta_{il}$ respectively in $(3.11b)$. Let $\bar{s}_i(x,y)$ and $\bar{s}_{il}(x,y)$ respectively denote $\bar{s}(x,y)$ in the transformed coordinates.
 In the former case we get a term
$$e^{-\lambda \cdot X_0}\int_{\{(x,y) \in D_i':\, |H_i(x,y)| < |\lambda|^{-{1 \over 100}}{\,\, {\rm or}\,\,} |\eta_i(x,y)| < |\lambda|^{-{1 \over 100}}\}} e^{-i\lambda_1 x \pm i\lambda_2y - \lambda_2k_i(x) - i\lambda_3s_i(x,y)} $$
$$\times |\bar{s}_i(x,y)|^z |H_i(x,y)|^{\delta z} e^{z^2} \rho_i(x,y)\,dx\,dy \eqno (3.12a)$$
Here $k_i(x)$ is as in the statement of Theorem 2.1, $H_i = H \circ \eta_i$, $\rho_i = \rho \circ \eta_i$, and so on. In the case of the $D_{il}'$ we 
analogously write
$$ e^{-i\lambda \cdot X_0}\int_{\{(x,y) \in D_{il}':\, |H_{il}(x,y)| < |\lambda|^{-{1 \over 100}}{\,\, {\rm or}\,\,} |\eta_{il}(x,y)| < |\lambda|^{-{1 \over 100}}\}} e^{-i\lambda_1 x \pm i\lambda_2y - \lambda_2k_{il}(x) - i\lambda_3s_{il}(x,y)} $$
$$\times |\bar{s}_{il}(x,y)|^z |H_{il}(x,y)|^{\delta z} e^{z^2} \rho_{il}(x,y)\,dx\,dy \eqno (3.12b)$$
Without loss of generality in our arguments we may take the $\pm i\lambda_2$ in $(3.12a)-(3.12b)$ to be $-i\lambda_2$. It will be better for our arguments 
if in the linear term in $k_i(x) = p_i x + l_i x^{s_i} + o(x^{s_i})$ is combined with the $x$ in the $\lambda_1 x$ term in $(3.12a)$, with the analogous
statement for $(3.12b)$. Thus we write $\lambda_0 = \lambda_1 + p_i\lambda_2$ and $K_i(x) = k_i(x) - p_ix$. So $K_i(x)$ has a zero of order
greater than $1$ at $x = 0$. Then $(3.12a)$ becomes
$$ e^{-i\lambda \cdot X_0}\int_{\{(x,y) \in D_i':\, |H_i(x,y)|< |\lambda|^{-{1 \over 100}}{\,\, {\rm or}\,\,}|\eta_i(x,y)| < |\lambda|^{-{1 \over 100}}\}} e^{-i\lambda_0 x - i\lambda_2y - \lambda_2K_i(x) - i\lambda_3s_i(x,y)} $$
$$\times |\bar{s}_i(x,y)|^z |H_i(x,y)|^{\delta z} e^{z^2} \rho_i(x,y)\,dx\,dy \eqno (3.13a)$$
We get an an analogous expression for the integrals $(3.12b)$, which we write as
$$ e^{-i\lambda \cdot X_0}\int_{\{(x,y) \in D_{il}':\,|H_{il}(x,y)| < |\lambda|^{-{1 \over 100}}{\,\, {\rm or}\,\,}|\eta_{il}(x,y)| < |\lambda|^{-{1 \over 100}}\}} e^{-i\lambda_0 x - i\lambda_2y - \lambda_2K_{il}(x) - i\lambda_3s_{il}(x,y)} $$
$$\times |\bar{s}_{il}(x,y)|^z |H_{il}(x,y)|^{\delta z} e^{z^2} \rho_{il}(x,y)\,dx\,dy \eqno (3.13b)$$
Note that by our constructions in section 3, $K_i(x)$ will be identically zero in $(3.13a)$ when 
$\beta_i = 0$ in Theorem 2.1; the cases where $k_i(x)$ is not linear and $\beta_i = 0$ are exactly the situations where one does the second resolution of
singularities using Theorem 2.2 and obtains the regions $D_{il}'$.

The strategy then will be as follows.  In the case where $\beta_i \geq 2$ after the application of Theorem 2.1 to $s(x,y)$, which occurs only
for terms of the form $(3.13a)$, we do a dyadic decomposition in $(3.13a)$ in both the $x$ and 
$y$ variables and then estimate each piece separately. We will apply the Van der Corput Lemma 3.2 for second derivatives in the $y$ direction on each dyadic piece,
and add over all pieces. We get a decay rate of $C|\lambda|^{-{1 \over 2}}$ from the second derivative Van der Corput lemma, times an additional
$C|\lambda|^{-\nu}$ for a small $\nu > 0$ which enables us to apply Theorem 3.1. This  $C|\lambda|^{-\nu}$ factor can come from one or both of
two places. First, on the domain where $|(x,y)| \sim |\eta_i(x,y)| < |\lambda|^{-{1 \over 100}}$, because $Re(z)$ is strictly greater than 
$\max (0, {1 \over 2} - g)$ and because $s_i(0,0) = 0$, one gets the additional $C|\lambda|^{-\nu}$ due to the $|\bar{s}_i(x,y)|^z$ damping factor in $(3.13a)$. Secondly, on the 
domain where $|H_i(x,y)| < |\lambda|^{-{1 \over 100}}$, the additional $C|\lambda|^{-\nu}$ factor arises from the $|H_i(x,y)|^{\delta z}$ factor in
$(3.13a)$. 

In the case where $\beta_i = 1$ after applying Theorem 2.1 to $s(x,y)$, which again only occurs in terms of the form $(3.13a)$, one applies the mixed derivative Van der Corput lemma 3.3 and one argues much like in the $\beta_i \geq 2$ situation.  If $\beta_i = 0$ after applying Theorem 2.1 to $s(x,y)$,
then there are two possibilities. First, $K_i(x)$ can be equal to zero in $(3.13a)$. In this case one does a dyadic decomposition in $x$ only, then 
applies Lemma 3.2 for second derivatives in the
$x$ direction instead of $y$, and argues as in the $\beta_i \geq 2$ case.

 If $\beta_i = 0$ after applying Theorem 2.1 to $s(x,y)$, but $K_i(x)$ is not identically zero, then our term is necessarily of the form $(3.13b)$. Things are 
more difficult here since a single application of a second derivative Van der Corput lemma will not suffice on some dyadic pieces. However, one
can do the following. To simplify the discussion, we assume the lower boundary of $D_{il}'$ is the $x$-axis. (The argument for the more general situation is
not fundamentally different). Letting $P_{\lambda}(x,y)$ denote the phase function $\lambda_0 x + i\lambda_2y + \lambda_2K_{il}(x) + \lambda_3s_{il}(x,y)$, 
observe that
$$\partial_{xx} P_{\lambda} (x,y) = \lambda_2 K_{il}''(x) + \lambda_3 \partial_{xx} s_{il}(x,y) \eqno (3.14)$$
By Theorem 2.1, one has that
$\partial_{xx} s_{il}(x,y) \sim x^{\alpha_i - 2}$ and $\partial_{xxx} s_{il}(x,y) \sim x^{\alpha_i - 3}$ ($\alpha_i \neq 2$ in these situations.) In addition,
the leading term of the Taylor series of $K_{il}(x)$ is of the form $l_{il} x^{s_i}$ for some $s_i < \alpha_i$. (There is no $l$ dependence in the exponent
$s_i$ since the exponent $s_i$ does not change after
the application of Theorem 2.2.)
Thus $\partial_{xx} P_{\lambda} (x,0) = \lambda_2 K_{il}''(x) + \lambda_3 \partial_{xx} s_{il}(x,0)$ has at
most one zero in $x$. Denote this value of $x$, if it exists, by $\tilde{x}$. There will be a certain $0 < r < 1$ such that if $x/\tilde{x} < r$ or
$x/\tilde{x}  > {1 \over r}$, then 
$\partial_{xx} P_{\lambda} (x,y)  > C |\lambda| x^{\alpha_i - 2}$, in which case one can use the Van der Corput lemma for second $x$ derivatives, similarly 
to how one did in the above case where $\beta_i = 0$ and $K_i(x)$ is identically zero.

When $r < \tilde{x}/x  < {1 \over r}$, one does a dyadic decomposition in $x$ and $y$ over pieces of the form $|x - \tilde{x}| \sim 2^{-j}, y \sim
2^{-k}$. Carefully examining the Taylor expansion of $s_{il}(x,y)$, we will see that on each such dyadic piece, we will be able to use a Van der Corput style lemma for second derivatives, either in the 
$x$ direction, the $y$ direction, or for a mixed second partial, such that adding over all of these dyadic pieces gives the bound of
 $C|\lambda|^{-{1 \over 2} - \nu}$ needed to apply Theorem 3.1.
The analysis will draw on Lemmas 2.5 and 2.6 and is arguably the most technically difficult segment of the proof of Theorem 1.1. Roughly speaking,
the phase function in this situation will effectively be of the form
$$\lambda_0x + \lambda_2y + \lambda_3c_1x^{\alpha_i - 3}(x - \tilde{x})^3 + \lambda_3c_2x^{p_{il}}y + \lambda_3c_3x^{\gamma_{il}}y^{\delta_{il} + 2}$$
Here $p_{il} \geq 1$ and $\gamma_{il},\delta_{il} \geq 0$ where the lowest edge of the Newton polygon of $s_{il}(x,y)$ connects $(\alpha_i,0)$
to $(\gamma_{il},\delta_{il} + 2)$. Because we are using Van der Corput lemmas for second derivatives, the $\lambda_0x + \lambda_2y$ plays no role here. The combined effect of the remaining three terms will enable us to use the Van der Corput lemmas in the desired fashion.

The above argument doesn't work when the upper vertex of the lowest non-horizontal edge of the Newton 
polygon of $s_{il}(x,y)$ is of the form $(A_{il}, 1)$, mainly because Lemma 2.5 doesn't apply to this situation. However, we will see that this case can be
dealt with using the mixed derivative Van der Corput Lemma 3.3, quite similarly to the $\beta_i = 1$ situation described above. 

The sharpness statement given by part b) of Theorem 1.1 will be shown via an explicit example.

\subsection {Overview of the proof of Theorem 1.2.}

The proof of Theorem 1.2 will essentially be a somewhat simpler version of the proof of Theorem 1.1. Because Radon transforms are translation invariant, it
suffices to assume the surface $S$ in $(1.2)$ is centered at the origin. In other words, we may assume that $(x_0,y_0,z_0) = (0,0,0)$ in $(1.5)$. Then
$(1.5)$ becomes
$$Rf(x,y,z) = \int_{\R^2} f(x - x',y - y', z - s(x',y')) \phi(x',y')\,dx' dy' \eqno (3.15) $$
This time, we embed $R$ in the analytic family $R_z$, where
$$R_zf(x,y,z) =\int_{\R^2} f(x - x',y - y', z - s(x',y')) |\bar{s}(x,y)|^z  e^{z^2} \phi(x',y')\,dx' dy' \eqno (3.16) $$
Note that $R_0 = R$. We will see that if $s > -g$ is fixed, then $||R_z f||_{L^p} \leq C||f||_{L^p}$ for arbitrarily large finite $p$, where 
$C$ is uniform over fixed $s$. Hence $R$ gains at least zero $L^2$ Sobolev derivatives on $L^p$. Analogous to the case of the maximal averages, this will
be proven simply by bounding the integral of $|\bar{s}(x,y)|^z$ and observing that an added $|e^{z^2}|$ factor makes the bound uniform over $Re\,z = s$. 
We don't use $L^{\infty}$ here since complex interpolation does not work well with $L^{\infty}$ Sobolev spaces.

We will then show for $s > \max(0, {1 \over 2} - g)$ and $Re\,z = s$ that one has the estimate $||R_z f||_{L^2_{1 \over 2}} \leq C||f||_{L^2}$
 with a constant $C$ that is uniform over $Re\,z = s$. This will be done by looking at the Fourier transform of the surface measure in $(3.16)$, 
given by 
$$\int_{\R^2} e^{-i\lambda_1 x - i\lambda_2y - i\lambda_3s(x,y)} |\bar{s}(x,y)|^z  e^{z^2}  \phi(x,y)\,dx\,dy \eqno (3.17)$$
This is the same as $(3.9')$, other than the removal of a magnitude 1 factor in front, a different-named cut-off function, and most importantly, with the
$|H(x,y)|^{-\delta z}$ removed. The effect of the removal of the $|H(x,y)|^{-\delta z}$ factor is that in imitating the above analysis following $(3.9')$, 
we will get a bound of $C_s|\lambda|^{-{1 \over 2}}$ instead of  $C_s|\lambda|^{-{1 \over 2} - \nu_s}$. Since we are looking to gain exactly ${1 \over 2}$
derivatives, this will give the desired estimates. 

Using complex interpolation between the two vertical lines above will then give Theorem 1.2a) for $p > 2$. Using duality about $p = 2$ will then give it for
remaining $p$. The sharpness statement of part b) will follow relatively quickly from a sharpness statement from [G6].

\section  {The proofs of Theorems 1.1 and 1.2.}

\subsection {The proof of Theorem 1.1.}

\subsubsection{\bf $L^{\infty}$ to $L^{\infty}$ boundedness for $Re\,z > -g$.}

\noindent Using the notation of $(3.8)$, for $z = s + it$ with $s > - g$ we first define $E(z)$ by
$$E(z) =\int_{\R^2} |\bar{s}({\bf x} - {\bf x_0})|^s|H({\bf x} - {\bf x_0})|^{\delta s}  e^{s^2 - t^2} |\phi(x,y)|\,dx\,dy \eqno (4.1)$$
$E(z)$ is the $L^1$ norm of the density function of the maximal function $M_z f$ as written in $(3.8)$. 
Shifting variables by ${\bf x_0}$ in $(4.1)$ and letting $r$ be such that 
$\rho(x,y)$ is supported on $|(x,y)| < r$ we get
$$E(z) =\int_{\{(x,y): |(x,y)| < r\}} |\bar{s}(x,y)|^s|H(x,y)|^{\delta s}  e^{s^2-t^2} |\rho(x,y)|\,dx\,dy \eqno (4.1')$$
 Using  $(3.10)$ we  get
$$E(z) \leq A e^{s^2 - t^2} \int_{\{(x,y): |(x,y)| < r\}} |\bar{s}(x,y)|^s|H(x,y)|^{\delta s}|(x,y)|^{-\beta} \,dx\,dy \eqno (4.2)$$
Suppose $0 > s > -g$ and $\delta$ is small enough that $|H(x,y)|^{\delta s}|(x,y)|^{-\beta}$ is integrable over $\{(x,y): |(x,y)| < r\}$. Suppose $p > 1$ is such that
$|H(x,y)|^{\delta p s}|(x,y)|^{-\beta}$ is also integrable over $\{(x,y): |(x,y)| < r\}$, and let $p'$ satisfy ${1 \over p} + {1 \over p'} = 1$. Define the
constant $C_{\delta,p,s}$ by
$$C_{\delta,p,s} = \bigg(\int_{\{(x,y): |(x,y)| < r\}} |H(x,y)|^{\delta p s}|(x,y)|^{-\beta}\,dx\,dy \bigg)^{1 \over p} \eqno (4.3)$$
 Then by H\"{o}lder's inequality applied with the measure 
$|(x,y)|^{-\beta}\,dx\,dy$, we have
$$E(z) \leq C_{\delta,p,s} \, e^{s^2 - t^2} \bigg(\int_{\{(x,y): |(x,y)| < r\}} |\bar{s}(x,y)|^{p's}|(x,y)|^{-\beta} \,dx\,dy \bigg)^{1 \over p'}\eqno (4.4)$$
We write the integral in $(4.4)$ as
$$\sum_{j = 0}^{\infty} \int_{\{(x,y): |(x,y)| < r,\, 2^{-j-1} \leq |\bar{s}(x,y)| < 2^{-j}\}} |\bar{s}(x,y)|^{p's}|(x,y)|^{-\beta} \,dx\,dy \eqno (4.5)$$
Using the definition $(1.6)$ of $g$ when $\bar{s}(x,y) = s(x,y)$, and Lemma 2.6 otherwise, we have that  the quantity $(4.5)$ is bounded by
$$C \sum_{j = 0}^{\infty} 2^{-jp's} \times j 2^{-jg} \eqno (4.6)$$
Using that $s < 0$, we see that as long as $p' < -{g \over s}$, the sum $(4.6)$ is finite. As a result $E(z)$ is uniformly bounded over $z$ with $Re\,z = s$.
 Because $E(z)$ is the $L^1$ norm of the density function of the maximal function $M_z f$, this means that if $p' < -{g \over s}$ then one has 
estimates $||M f||_{\infty} \leq C||f||_{\infty}$ with constant $C$ uniform over $z$ with $Re\,z = s$. Whenever $0 > s > -g$, we have that $-{g \over s} > 1$
and there will be some value of $p$ for which $1 < p' <  -{g \over s} $. For this value of $p$ one can choose $\delta$ small enough for the above argument to 
work. As a result, whenever $0 > s > - g$ we have $||M f||_{\infty} \leq C||f||_{\infty}$ with constant $C$ uniform over $z$ with $Re\,z = s$, as needed.

\subsubsection {Fourier transform estimates when $|\lambda_1| + |\lambda_2| > |\lambda_3|$.}

We assume $z$ is such that $s = Re\,z > \max(0, {1 \over 2} - g)$, and we will bound the surface measure Fourier transform $(3.9')$ under the assumption that 
$|\lambda_1| + |\lambda_2| > |\lambda_3|$. Note that $g$ is maximized when $s(x,y)$ is the nondegenerate function $x^2 + y^2$, and a simple
calculation reveals $g = 1 - {\beta \over 2}$ in this case. Thus we always have $s > {1 \over 2} - g \geq {1 \over 2} - (1 - {\beta \over 2}) = {\beta - 1 \over 2}$. Since $\bar{s}(x,y)$ has a zero of order at least $2$ 
at the origin, this means that if $s > \max(0, {1 \over 2} - g)$ then we have 
$$|\bar{s}(x,y)|^s \leq C|(x,y)|^{2s} \leq C|(x,y)|^{\beta - 1} \eqno (4.7)$$
Since $|H(x,y)|^{\delta s}$ is bounded  and $|\rho(x,y)| \leq A|(x,y)|^{-\beta}$ by $(3.10)$, if $s > \max(0, {1 \over 2} - g)$ then we therefore have
$$|\bar{s}(x,y)|^s |H(x,y)|^{\delta s} |\rho(x,y)| \leq C'|(x,y)|^{-1} \eqno (4.8)$$
We break the integral $(3.9')$ into $|(x,y)| \leq {1 \over |\lambda|}$ and $|(x,y)| > {1 \over |\lambda|}$ parts. To estimate the first part, we simply
take absolute values of the integrand and integrate, using $(4.8)$. The result is a bound of $C|\lambda|^{-1}$. Thus we devote our attention to the second
term, which we denote by $J(\lambda)$. Thus we have
$$|J(\lambda)| =  \bigg| \int_{|(x,y)| > |\lambda|^{-1}} e^{-i\lambda_1 x - i\lambda_2y - i\lambda_3s(x,y)} |\bar{s}(x,y)|^z |H(x,y)|^{\delta z} e^{z^2}\rho(x,y)\,dx\,dy
\,\bigg|\eqno (4.9)$$
We write $J(\lambda) =\sum_{i=1}^4 J_i(\lambda)$, where $J_i(\lambda)$ is the portion of $(4.9)$ in one of the four quadrants.
Assuming we are in a sufficiently small neighborhood of the origin, there is a finite list of directions $\{v_j\}_{j=1}^N$ not in the $x$ or $y$ direction
 and a constant $C$ such that on the domain of integration of each $J_i(\lambda)$ there is some $v_i$ such that 
$$|\partial_{v_i} (\lambda_1 x + \lambda_2y + \lambda_3s(x,y))| > C|\lambda| \eqno (4.10)$$
One can choose $v_i$ such that $(4.10)$ holds due to the condition that $|\lambda_1| + |\lambda_2| > |\lambda_3|$ and the fact that $\nabla s(0,0) = 0$.

We would like to now apply the Van der Corput lemma, Lemma 3.2, for first derivatives in the $v_i$ direction in $(4.9)$. Because $\bar{s}(x,y)$ is built from finite-type functions, we can choose the $v_i$ such that each cross-section of the domain of integration in $(4.9)$ in the $v_i$ direction can be written as the union of boundedly many intervals on which $\bar{s}(x,y), \partial_{v_i}\bar{s}(x,y)$, and $\partial_{v_iv_i}\bar{s}(x,y)$ are nonvanishing. If $H(0,0) = 0$, we may also assume each such cross-section is the union of boundedly many intervals on which $H(x,y)$ and $\partial_{v_i}H(x,y)$ 
are nonvanishing. As a result, each cross-section is the union of boundedly many intervals on which Lemma 3.2 for first derivatives applies and one can add
 the resulting estimates. 

To understand the estimate one obtains from this application of Lemma 3.2, one must examine the result of integrating the absolute value of 
 $\partial_{v_i} (|\bar{s}(x,y)|^z |H(x,y)|^{\delta z} \rho(x,y))$ in the $v_i$ direction. If the derivative lands on $\rho(x,y)$, one incurs an additional factor of $C|(x,y)|^{-1}$ due to $(3.10)$.
Using $(4.8)$ and the fact that the intervals of integration have length bounded by $C|(x,y)|$, the integral of this term is therefore bounded by $C|(x,y)|^{-1}$.

 If the derivative lands on $|\bar{s}(x,y)|^z$, the $|\bar{s}(x,y)|^z$ becomes
$z|\bar{s}(x,y)|^{z-1}\partial_{v_i}|\bar{s}(x,y)|$. Observe that by $(3.10)$ and the fact that $H(x,y)$ is bounded, we have
$$\big|z \bar{s}(x,y)|^{z-1}\partial_{v_i}|\bar{s}(x,y)\big| \times |H(x,y)|^{\delta z} |\rho(x,y)| \leq |z|| \bar{s}(x,y)|^{s-1}\partial_{v_i}|\bar{s}(x,y)|
\times |(x,y)|^{-\beta}$$
Since $|(x,y)|$ is within a constant factor of some fixed $|(x^*,y^*)|$ on any interval of integration, on such an interval we then have 
$$|z ||\bar{s}(x,y)|^{z-1}\partial_{v_i}|\bar{s}(x,y)|\times |H(x,y)|^{\delta z} |\rho(x,y)| \leq C |z ||\bar{s}(x,y)|^{s-1}\partial_{v_i}|\bar{s}(x,y)|\times |(x^*,y^*)|^{-\beta}\eqno (4.11)$$
The right-hand side of $(4.11)$ is a constant times $\big|\partial_{v_i} |\bar{s}(x,y)|^s\big| $ on any interval of integration in the $v_i$ direction in $(4.9)$. Hence on each 
such interval
we may integrate the derivative back to the original function at the endpoints, and what we get is $C {|z| \over s} |\bar{s}(x,y)|^s |(x^*,y^*)|^{-\beta}$
at the endpoints.  Since
we are assuming $s = Re\,z$ is fixed here, we can incorporate it into the constant $C$. Using $(4.7)$ this is at most $C|z| |(x^*,y^*)|^{-1}$.

If $H(0,0)$ and the derivative lands on the function $|H(x,y)|^{\delta z}$, similarly to the above situation one gets a term whose integral is bounded by $C|z| |(x^*,y^*)|^{-1}$. If $H(0,0) \neq 0$, the derivative simply causes a $C|z|$ factor to be incurred and the integral of the resulting term is bounded by $C|z|$.

Thus regardless of where the derivative lands, the absolute value of the resulting term integrates to $C|z| |(x^*,y^*)|^{-1}$. Thus when applying Lemma 
3.2 in the $v_i$ direction, this expression can be used as a bound for the term called $\int_a^b |\psi'(x)|\,dx$ in $(3.4)$. By $(4.8)$ it can also be used
for a bound for the $|\psi(b)|$ term there. Thus if we apply Lemma 3.2 in the $v_i$ direction in $(4.9)$ and then integrate the result in the $x$  
direction, one obtains
$$|J_i(\lambda)| \leq C''|ze^{z^2}| \bigg|\int_{|\lambda|^{-1}}^{\infty} {1 \over |\lambda|x} \,dx\bigg|$$
Performing the integration, and observing that $|ze^{z^2}|$ is bounded on any vertical line in the complex plane, we obtain
$$|J_i(\lambda)| \leq  C'''{\ln |\lambda| \over |\lambda|} \eqno (4.12)$$
 Here $C'''$ is independent of $Im\,z$ for fixed $s$. Adding $(4.12)$ over the four quadrants and adding the result to the $C|\lambda|^{-1}$ bound for the 
$|(x,y)| \leq {1 \over |\lambda|}$ portion, we see that if $|\lambda_1| + |\lambda_2| > |\lambda_3|$, then in the notation of $(3.9')$ we have shown that
$$|\widehat{d\sigma_z}(\lambda)| \leq C {\ln |\lambda| \over |\lambda|} \eqno (4.13)$$
This is stronger than the exponent $-{1 \over 2} -\nu$ needed to apply Theorem 3.1.

\subsubsection {Fourier transform estimates when $|\lambda_1| + |\lambda_2| \leq |\lambda_3|$, $|H(x,y)| > |\lambda|^{-{1 \over 100}}$, $|(x,y)|  > |\lambda|^{-{1 \over 100}}$.}

We now estimate $|I_1(\lambda)|$, where $I_1(\lambda)$ is as in $(3.11a)$. The argument is based on a similar argument in [G7]. The idea is as follows. 
In the case of nonvanishing Hessian determinant, one can get the traditional estimate 
 $|I_1(\lambda)| < C|\lambda|^{-1}$. On the support of the integrand
of $I_1(\lambda)$, the Hessian determinant $H(x,y)$ is at least $|\lambda|^{-{1 \over 100}}$, and we will see
that after an argument elaborating that of the case of nonvanishing Hessian determinant, 
we still get an estimate $|I_1(\lambda)| < C|\lambda|^{-{3 \over 5}}$, an estimate better than what is needed.

Since $s(x,y)$ is of finite type and we are proving a local result, shrinking our neighborhood of the origin if necessary
 we may let $u$ and $v$ be nonparallel directions
such that for some $k \geq 2$, $\partial_u^k s(x,y)$, $\partial_v^k s(x,y)$, $\partial_u^k \bar{s}(x,y)$, $\partial_v^k \bar{s}(x,y)$, and $\partial_u\partial_v^{k-1}s(x,y)$  are nonvanishing on some disk of radius $r$ centered at
the origin that contains the support of the integrand of $I_1(\lambda)$.
Similarly, if $H(0,0) = 0$, we may further assume that for some $k' \geq 1$ we have that $\partial_u^{k'} H(x,y)$ and $\partial_v^{k'} H(x,y)$
are nonvanishing.

 Let $a_1$ and $a_2$ be the constants defined by  $a_1 = 
-\partial_u ({\lambda_1\over \lambda_3}x + {\lambda_2\over \lambda_3}y)$ and $a_2 = -\partial_v 
({\lambda_1\over \lambda_3}x + {\lambda_2\over \lambda_3}y)$. Let $D_r$ be the disk of radius $r$
centered at the origin, where $r$ is as above. Define the
sets $A_1$, $A_2$, and $A_3$ by
$$A_1 = \{(x,y) \in D_r: |\partial_us(x,y) - a_1| > |\lambda|^{-{1 \over 3}}\} \eqno (4.14a)$$
$$A_2 = \{(x,y) \in D_r: |\partial_us(x,y) - a_1| \leq |\lambda|^{-{1 \over 3}},\,|\partial_vs(x,y) - a_2| >
|\lambda|^{-{1 \over 3}} \} \eqno (4.14b)$$
$$A_3 = \{(x,y) \in D_r: |\partial_us(x,y) - a_1| \leq |\lambda|^{-{1 \over 3}},\, |\partial_vs(x,y) - a_2| \leq
|\lambda|^{-{1 \over 3}} \} \eqno (4.14c)$$
Correspondingly, write the contributions to $I_1(\lambda)$ from $A_1$, $A_2$, and $A_3$ as $J_1(\lambda)$, $J_2(\lambda)$, and
$J_3(\lambda)$ respectively. To analyze $J_1(\lambda)$, we apply Lemma 3.2 for first derivatives in the $u$ direction. Since $\partial_u^k 
s(x,y)$ is nonvanishing, the condition $(4.14a)$ will cause there to be boundedly many intervals on which to apply the lemma. 

Note that $|H(x,y)|^{\delta z}$ and $|\bar{s}(x,y)|^z$ are bounded since $Re\,z > 0$ here and that $|\rho(x,y)| < C|\lambda|^{2 \over 100}$ on the domain of
integration due to $(3.10)$ and the condition that $|(x,y)| > |\lambda|^{-{1 \over 100}}$. As a result, the term called $|\psi(b)|$ in $(3.4)$ can
be taken to be $C|\lambda|^{2 \over 100}$.

Moving to the term called $\int_a^b |\psi'(x)|\,dx$ in $(3.4)$, we get several terms, depending on where the derivative lands. If it lands on the
$\rho(x,y)$ factor, we gain a factor of $C|(x,y)|^{-1} < C'|\lambda|^{1 \over 100}$ due to $(3.10)$, so that this term contributes 
$C|\lambda|^{3 \over 100}$ to $\int_a^b |\psi'(x)|\,dx$.

If  the derivative lands on the $|H(x,y)|^{\delta z}$ factor and $H(0,0) = 0$, we take absolute values and integrate in the $u$ direction, 
 using that $|\partial_u^{k'}H(x,y)|$ is
nonvanishing to ensure that there are boundedly many intervals on which $H(x,y)$ is nonzero and
monotone and thus on which we can integrate back its $u$-derivative. The result is a bound of $C|z||\lambda|^{2 \over 100}$.
Although we get a $C|z|$ factor here, the presence of the $e^{z^2}$ in the damping factor is more than enough to compensate. 
If $H(0,0) \neq 0$, then a derivative landing on $|H(x,y)|^{\delta z}$ has the effect of simply adding a $C|z|$ factor, which gives the same estimate 
as the $H(0,0) = 0$ situation.

Lastly, suppose the derivative lands on the factor $|\bar{s}(x,y)|^z$. The directions $u$ and $v$ were defined so that
$\bar{s}(x,y)$ has certain nonvanishing higher order derivatives in the $u$ and $v$ directions. So one can argue 
as above, breaking up the one-dimensional integration into boundedly many 
intervals on which $\bar{s}(x,y)$ is nonzero and monotone. Hence we get the same upper bounds as before. 

We have now considered all possible places the derivative can land, and we see that the $\int_a^b |\psi'(x)|\,dx$ term of $(3.4)$ is bounded by
$C|\lambda|^{3 \over 100}$. Applying Lemma 3.2 now, using the lower bounds on the $u$ derivative of the phase provided by $(4.14a)$, we see
that $|J_1(\lambda)| \leq C |\lambda|^{3 \over 100} \times |\lambda|^{-1} \times |\lambda|^{1 \over 3} < C|\lambda|^{-{3 \over 5}}$, the desired estimate.

The bounds for $|J_2(\lambda)|$ are proven exactly as they were for $|J_1(\lambda)|$, replacing the roles of the $u$ and $v$ 
variables. The presence of the added condition $|\partial_us(x,y) - a_1| \leq |\lambda|^{-{1 \over 3}}$ in 
the domain, which does not have an analogue above, does not interfere with any of the above estimates; 
the condition  that $\partial_u\partial_v^{k-1}s(x,y)$ is nonvanishing ensures that in
$v$ direction, one still has boundedly many intervals on which to apply Lemma 3.2.

We now move on to $J_3(\lambda)$. We write the domain of integration as the union of at most $C|\lambda|^{{2 \over 100}}$ subdomains,
each of which is the intersection of the original domain of integration with a square of diameter $c|\lambda|^{-{1 \over 100}}$, where $c$ is a constant
to be determined by our arguments.

Let $S$ be any of these squares. We consider the level sets of $\partial_u s(x,y)$ and $\partial_v s(x,y)$ on $S$. The gradients of 
both functions are bounded below in absolute value by $C|H(x,y)|$. Since in the case at hand $|H(x,y)| > |\lambda|^{-{1 \over 100}}$ on the 
domain of integration, one has  $|H(x,y)| > {1 \over 3}|\lambda|
^{-{1 \over 100}}$ on $S$ if we chose the constant $c$ in the diameter $c|\lambda|^{-{1 \over 100}}$ 
of the squares sufficiently small. As a result, if $c$ is small enough the level sets of both 
$\partial_u s(x,y)$ and $\partial_v s(x,y)$ do not self-intersect on $S$. Hence we may use $\partial_u s(x,y)$
and $\partial_v s(x,y)$ as coordinates on $S$. In particular, we may evaluate the measure of the set $A_3 \cap S$, where $A_3$ is as in 
of $(4.14c)$, by changing into these coordinates in the integral of its characteristic function. The result is
$$|A_3 \cap S| < C \min_S |H(x,y)|^{-1} |\lambda|^{-{2 \over 3}} \eqno (4.15)$$
So we conclude that $|A_3 \cap S| < C'|\lambda|^{-{2 \over 3} + {1 \over 100}}$. Adding over all $S$ we get $|A_3| <
 C'|\lambda|^{-{2 \over 3} + {3 \over 100}}$. Since the integrand of $J_3(\lambda)$ is uniformly
bounded on $Re(z) = s > 0$, we conclude that 
$$|J_3(\lambda)| \leq C''|\lambda|^{-{2 \over 3} + {3 \over 100}} < C''|\lambda|^{-{3 \over 5}} \eqno (4.16)$$
This gives the needed estimate. Adding the contributions from $|J_1(\lambda)|$, $|J_2(\lambda)|$, and $|J_3(\lambda)|$, we conclude
that  $|I_1(\lambda)|$ in $(3.11a)$ is at most $C'''|\lambda|^{-{3 \over 5}}$, for a constant $C'''$ independent of $Im\,z$ for fixed $Re\,z = s$.
Since $-{3 \over 5} < -{1 \over 2}$ we conclude that $|I_1(\lambda)|$ satisfies the bounds needed to apply Theorem 3.1.

\subsubsection {Fourier transform estimates when $|\lambda_1| + |\lambda_2| \leq |\lambda_3|$ and $|H(x,y)| < |\lambda|^{-{1 \over 100}}$ or  $|(x,y)|  < |\lambda|^{-{1 \over 100}}$.}

\noindent {\bf Case 1: $\beta_i \geq 2$.}

Here we are bounding $(3.13a)$, and in the case at hand $\bar{s}(x,y) = s(x,y)$. As always, we assume $Re\,z$ is some fixed $s > \max(0, {1 \over 2} - g)$. 
We split the integral $(3.13a)$ dyadically in both $x$ and $y$. Let $J_{jk}$ denote the interval
$[2^{-j-1}, 2^{-j}] \times [2^{-k-1}, 2^{-k}]$. So we are bounding $\sum_{j,k} |I_{jk}(\lambda)|$, where
$$ I_{jk}(\lambda) = \int_{\{(x,y) \in D_i' \cap J_{jk}:\, |H_i(x,y)|< |\lambda|^{-{1 \over 100}}{\,\, {\rm or}\,\,}|\eta_i(x,y)| < |\lambda|^{-{1 \over 100}}\}} e^{-i\lambda_0 x - i\lambda_2y - \lambda_2K_i(x) - i\lambda_3s_i(x,y)} $$
$$\times |s_i(x,y)|^z |H_i(x,y)|^{\delta z} e^{z^2} \rho_i(x,y)\,dx\,dy \eqno (4.17)$$
We will apply the Van der Corput Lemma 3.2 in the $y$ direction for second derivatives. Note that the second $y$ derivative of the phase function  in $(4.17)$ is
given by $-\lambda_3 \partial_{yy} s_i(x,y)$, and by $(2.1)$ one has $|\partial_{yy}s_i(x,y)| > Cx^{\alpha_i}y^{\beta_i - 2}$ on $D_i'$. This is the
second derivative lower bound that we will use in Lemma 3.2.

The rotation performed at the beginning of section 2 ensures that $\partial_y^o s(x,y)$ is nonzero, where $o$ is the order of the zero of $s(x,y)$ at the origin, and that 
$\partial_y^p H_i(x,y)$ is nonzero for some $p \geq 1$ if $H_i(0,0) = 0$. As a result, for fixed $x$ there are boundedly 
many intervals in the $y$ direction on which $s_i(x,y)$ and  $\partial_y s_i(x,y)$ are nonzero, and the same is true for $H_i(x,y)$ and $\partial_y H_i(x,y)$
if $H_i(0,0) = 0$.
 We apply Lemma 3.2 on each of these  intervals
and add the results. To this end, we must bound the quantity $(3.4)$. First, writing $z = s + it$, note by $(3.10)$ that 
$$\big||s_i(x,y)|^z |H_i(x,y)|^{\delta z} e^{z^2}\rho_i(x,y)\big| \leq A|e^{z^2}| |s_i(x,y)|^s |H_i(x,y)|^{\delta s}|(x,y)|^{-\beta} \eqno (4.18)$$
Next, we examine the $y$ derivative of $|s_i(x,y)|^z |H_i(x,y)|^{\delta z} e^{z^2}\rho_i(x,y)$. We get several terms, depending on which factor 
the derivative lands. Suppose it lands on the $|s_i(x,y)|^z$. Note that $\partial_y (|s_i(x,y)|^z) =  \pm z |s_i(x,y)|^{z-1}\partial_y s_i(x,y)$.
Since $|s_i(x,y)| \sim x^{\alpha_i}y^{\beta_i}$ and $|\partial_y s_i(x,y)|
 \sim x^{\alpha_i}y^{\beta_i-1}$ by $(2.2)$, one therefore has that 
$$\big|\partial_y (|s_i(x,y)|^z)\big| \leq C|z|{1 \over y}|s_i(x,y)|^z \eqno (4.19)$$
Suppose the derivative lands on $\rho_i(x,y)$. By $(3.10)$ and the fact that the coordinate change $\eta_i(x,y)$ is of the form $(x,y + k_i(x))$ where 
$k_i(x)$ has a zero of order at least one at the origin, $\rho_i(x,y)$ also satisfies $(3.10)$. Thus the derivative results in an additional factor of 
$C{1 \over x}$. Since $|y| < C'|x|$ on all $D_i$, this is better than incurring a factor of $C{1 \over y}$. Thus in view of $(3.10)$ we have
$$|\partial_y \rho_i(x,y)| \leq C{1 \over y}|(x,y)|^{-\beta - 1} \eqno (4.20)$$
Comparing $(4.19)$ and $(4.20)$ to $(4.18)$ in the context of $(3.4)$, we see that when seeking to apply Lemma 3.2 to the $y$ integral in the dyadic 
rectangle, the terms where the derivative lands on $\rho_i(x,y)$ or $|s_i(x,y)|^z$ give a contribution in the integral term of $(3.4)$ that is no worse
 than $C|z|$ times the bound given by $(4.18)$ for the term denoted by $\psi(b)$ in $(3.4)$. Since  $|s_i(x,y)| \leq C2^{-j\alpha_i}2^{-k\beta_i}$ and
$|\rho_i(x,y)| \leq C|(x,y)|^{-\beta} \leq C'2^{-j\beta}$ (recall $|y| < C|x|$ on all of our domains), we may write this bound in the form
 $$C |ze^{z^2}| \sup_{J_{jk}}|H_i(x,y)|^{\delta s} \times 2^{-js\alpha_i}2^{-ks\beta_i} \times  2^{-j\beta} \eqno (4.21)$$ 
Lastly, we consider the term where the derivative lands on  $|H_i(x,y)|^{\delta z}$. First suppose $H_i(0,0) = 0$. We bound all the remaining factors by constants as in $(4.21)$,
and then integrate the resulting function $C''\big|\partial_y |H_i(x,y)|^{\delta z} \times
2^{-j\alpha_i}2^{-k\beta_i} \times  2^{-j\beta}\big|$ back to $C'''|z| |H_i(x,y)|^{\delta s} \times 2^{-js\alpha_i}2^{-ks\beta_i} \times  2^{-j\beta}$, to once again get a bound of $(4.21)$. In the case where $H_i(0,0) \neq 0$ the derivative just causes an additional $C|z|$ factor to be incurred, so $(4.21)$ still holds.

Combining all of the above, we see that in the situation at hand, the expression in parentheses on the right in $(3.4)$ is bounded by $(4.21)$.
Since $ze^{z^2}$ is bounded on any vertical strip in the complex plane, there is a constant $C_s$ independent of $Im\, z$ such that we may write the 
bound as
$$C_s \sup_{J_{jk}}|H_i(x,y)|^{\delta s} \times 2^{-js\alpha_i}2^{-ks\beta_i} \times  2^{-j\beta} \eqno (4.21')$$
We now are in a position to apply Lemma 3.2 in the $y$ direction. Since $\beta_i > 1$, we have $|\partial_{yy} s_i(x,y)| > c 2^{-j\alpha_i}2^{-k(\beta_i - 2)}$
for some constant $c$. Note that since $|\lambda_1| + |\lambda_2| \leq |\lambda_3|$, we have $|\lambda| \sim |\lambda_3|$. Applying Lemma 3.2 
in the $y$ direction and integrating the result in $x$, we obtain
$$|I_{jk}| \leq C_s' {1 \over (|\lambda| 2^{-j\alpha_i}2^{-k(\beta_i - 2)})^{1 \over 2}} \sup_{J_{jk}}|H_i(x,y)|^{\delta s} \times 2^{-js\alpha_i}2^{-ks\beta_i} \times  2^{-j\beta} \times 2^{-j} \eqno (4.22)$$
It is more convenient for our arguments to write this in the form
$$|I_{jk}| \leq C_s''|\lambda|^{-{1 \over 2}}\sup_{J_{jk}}|H_i(x,y)|^{\delta s}\int_{J_{jk}}|s_i(x,y)|^{s - {1 \over 2}}|(x,y)|^{-\beta}\,dx\,dy \eqno (4.23a)$$
Since $s > \max(0, {1 \over 2} - g)$, we may write $s = {1 \over 2} - g^*$ where $g^* < \min({1 \over 2},g)$. In this notation, $(4.23a)$ becomes
$$|I_{jk}| \leq C_s''|\lambda|^{-{1 \over 2}}\sup_{J_{jk}}|H_i(x,y)|^{\delta s}\int_{J_{jk}}|s_i(x,y)|^{-g^*}|(x,y)|^{-\beta}\,dx\,dy \eqno (4.23b)$$
We next divide the $I_{jk}$ into two types, in accordance with the domain of integration of $(4.17)$. The first type of $I_{jk}$ are those for which 
$|(x,y)| >  C_0|\lambda|^{-{1 \over 100}}$ holds on the whole interval $J_{jk} = [2^{-j-1}, 2^{-j}] \times [2^{-k-1}, 2^{-k}]$, where $C_0$ is chosen such
that in this situation $|\eta_i(x,y)| > |\lambda|^{-{1 \over 100}}$ on $J_{jk}$; this can be done since by the form of $\eta_i(x,y)$ given in Theorem 2.1, we always have $|(x,y)| \sim |\eta_i(x,y)|$. Given the form of the integral $(4.17)$, on the domain of integration of any $I_{jk}$ of the first type one has $|H_i(x,y)| < |\lambda|^{-{1 \over 100}}$. The second type of $I_{jk}$ are simply those which are not of the first type. On the domain of integration of 
the second type of $I_{jk}$ one has $|(x,y)| < C_1|\lambda|^{-{1 \over 100}}$ for some constant $C_1$.

For $I_{jk}$ of the first type, we insert $|H_i(x,y)| < |\lambda|^{-{1 \over 100}}$ into $(4.23b)$ and add over all such $I_{jk}$. The result is a bound of
$$C_s''' |\lambda|^{-{1 \over 2} - {\delta \over 100}}\int_{|(x,y)| < r} |s(x,y)|^{-g^*}|(x,y)|^{-\beta}\,dx\,dy \eqno (4.24)$$
Here $r$ is small enough so that the support of $\rho_i(x,y)$ is contained in the domain of integration of $(4.24)$. The definition $(1.6)$ of $g$ implies
that if $r$ is sufficiently small, then $\int_{|(x,y)| < r} |s(x,y)|^{-h}|(x,y)|^{-\beta}\,dx\,dy$ is finite for $h < g$ and infinite for $h > g$. Since 
$g^* < g$, $(4.24)$ gives a bound of $C_s''' |\lambda|^{-{1 \over 2} - {\delta \over 100}}$. As this exponent is less than $-{1 \over 2}$, with a 
constant independent of $Im\,z$, this
suffices for our application of Theorem 3.1.

We now move to $I_{jk}$ of the second type. Let $g^{**}$ satisfy $g^* < g^{**} < \min({1 \over 2}, g)$. Since $|s(x,y)| < C|(x,y)|^2$ and 
$|(x,y)| < C_1|\lambda|^{-{1 \over 100}}$ on the domain of integration on an $I_{jk}$ of the second type, on this domain of integration we have
$$|s(x,y)|^{-g^*} = |s(x,y)|^{g^{**} - g^*} |s(x,y)|^{-g^{**}} \leq C|\lambda|^{-{g^{**} - g^{*} \over 50}} |s(x,y)|^{-g^{**}} \eqno (4.25)$$
We insert this into $(4.23b)$, and use the fact that $H_i(x,y)$ is a bounded function to conclude that
$$|I_{jk}| \leq C_s |\lambda|^{-{1 \over 2}-{g^{**} - g^{*} \over 50}} \int_{I_{jk}}|s_i(x,y)|^{-g^{**}}|(x,y)|^{-\beta}\,dx\,dy \eqno (4.26)$$
We add this over all $I_{jk}$ of the second type. Using that $g^{**} < g$, we get a bound of
$C_s' |\lambda|^{-{1 \over 2}-{g^{**} - g^{*} \over 50}}$
Since the exponent here is once again less than $-{1 \over 2}$, with a constant independent of $Im\,z$, this again suffices to apply Theorem 3.1.
Adding the estimates obtained above over $I_{jk}$ of both types gives a bound of the form $C_s|\lambda|^{-{1 \over 2} - \nu_s}$ for $(3.13a)$ as needed.

\noindent {\bf Case 2: $\beta_i = 1$.}

Like above, $\bar{s}(x,y) = s(x,y)$ here. We once again decompose $(3.13a)$ dyadically in both $x$ and $y$, and  bound $\sum_{j,k} |I_{jk}(\lambda)|$, where
$I_{jk}(\lambda)$ is as in $(4.17)$. This time, we use Lemma 3.3 in place of Lemma 3.2. Because the damping function in $(4.17)$ is not necessarily a $C^1$ function of
$y$ for fixed $x$, strictly speaking Lemma 3.3 does not immediately apply, but the proof of Lemma 3.3 in [G1] works equally well  if for fixed $x$ the damping function is 
a piecewise $C^1$ function where the number of pieces is uniformly bounded. This is indeed the situation at hand.

By $(2.1)$ we have $s_i(x,y) \sim 2^{-j \alpha_i}2^{-k}$ and $\partial_{xy} s_i(x,y) \sim 2^{-j(\alpha_i - 1)}$ on the domain  of integration of $(4.17)$. 
As a result, if $P_{\lambda}(x,y)$ denotes the phase function in $(4.17)$ then we have $|\partial_{xy} P_{\lambda}(x,y)| \sim |\lambda| 2^{-j(\alpha_i - 1)}$ on the domain  of integration. 

The quantity called 
$\int_{\{y: (x,y) \in R\}} |\partial_y\Psi(x,y)|\,dy$ in $(3.5b)$ is exactly the integral in $y$ of $\partial_y(|s_i(x,y)|^z |H_i(x,y)|^{\delta z} e^{z^2}\rho_i(x,y))$ that was
computed above. Thus the quantity $N$ of $(3.5b)$ satisfies the same bound that the quantity $|\psi(b)| + \int_a^b |\psi'(x)|\,dx$ of $(3.4)$ was computed to satisfy above.
Namely, we have that $N$ is bounded by the right-hand side of $(4.21')$. Applying Lemma 3.3 now, we get that
$$|I_{jk}(\lambda)| \leq C_s {1 \over (|\lambda| 2^{-j(\alpha_i - 1)})^{1 \over 2}} \sup_{J_{jk}}|H_i(x,y)|^{\delta s} \times 2^{-js\alpha_i}2^{-s\beta_i} \times  2^{-j\beta} \times 
2^{-{j \over 2}} \times 2^{-{k \over 2}} \eqno (4.27)$$
This reduces to exactly $(4.23a)$. Arguing exactly as in the steps following $(4.23a)$ gives a bound of $C_s|\lambda|^{-{1 \over 2} - \nu_s}$ for $(3.13a)$ once again.

\noindent {\bf Case 3: $\beta_i = 0$ and $\eta_i$ is linear.}

In this case, $K_i(x) = 0$ in $(3.13a)$. As a result, the second $x$ derivative of the phase in $(3.13a)$ is the given by that of the 
$\lambda_3 s_i(x,y)$ term, and thus by $(2.1)$ has magnitude $\sim |\lambda_3x^{\alpha_i - 2}| \sim |\lambda x^{\alpha_i - 2}|$. We decompose
 $(3.13a)$ dyadically in the $x$ variable only this time. Let $L_j$ be the interval $[2^{-j-1},2^j]$. We will bound $\sum_j|N_j(\lambda)|$,
  where  $N_j(\lambda)$ is given by
  $$\int_{\{(x,y) \in D_i' \cap L_j :\, |H_i(x,y)|< |\lambda|^{-{1 \over 100}}{\,\, {\rm or}\,\,}|\eta_i(x,y)| < |\lambda|^{-{1 \over 100}}\}} e^{-i\lambda_0 x - i\lambda_2y -  i\lambda_3s_i(x,y)} $$
$$\times |s_i(x,y)|^z |H_i(x,y)|^{\delta z} e^{z^2} \rho_i(x,y)\,dx\,dy \eqno (4.28)$$
We apply Lemma 3.2 for second derivatives in the $x$ direction in $(4.28)$.  Although we are integrating in the $x$ direction, we can bound the
quantity $(3.4)$ in much the same way as we did for $y$ derivatives in $(4.18)-(4.21')$. For $(3.10)$ implies that our bounds for
$|\partial_x \rho_i(x,y)|$ are the same as our bounds for $|\partial_y \rho_i(x,y)|$. Taking an $x$ derivative of $|s_i(x,y)|^z$ induces a $C|z| {1 \over x}$
factor analogous to the $C|z| {1 \over y}$ factor in $(4.19)$, and an $x$ derivative landing on $|H_i(x,y)|$ can be dealt with exactly as a $y$ derivative landing on $|H_i(x,y)|$ was dealt with in the paragraph above $(4.21')$. In particular, the real-analyticity of $H_i(x,y)$ in $x^{1 \over N}$ and $y$ will ensure we always have boundedly many intervals of integration when $H_i(0,0) = 0$. The result is that in the case at hand, the expression in parentheses on the right of 
$(3.4)$ is given by
$$C_s \sup_{x \in L_j}|z||H_i(x,y)|^{\delta s} \times 2^{-js\alpha_i}\times  2^{-j\beta} \eqno (4.29)$$
Applying Lemma 3.2 in the $x$ direction in $(4.28)$ and then integrating the result in $y$, using that the second dervative of the phase is $\sim |\lambda x^{\alpha_i - 2}|$, gives
$$|N_j(\lambda)| \leq C_s' {1 \over (|\lambda| 2^{-j(\alpha_i - 2)})^{1 \over 2}} \sup_{x \in L_j}|H_i(x,y)|^{\delta s} \times 2^{-js\alpha_i} \times  2^{-j\beta} \times 2^{-jM_i} \eqno (4.30)$$
Here $M_i$ is as in Theorem 2.1; that is, the cross-sectional width of $D_i$ for fixed $x$ is $\sim x^{M_i}$ when $\beta_i = 0$. This time, it is more convenient 
to write this in the form
$$|N_j(\lambda)| \leq C_s'''|\lambda|^{-{1 \over 2}}\sup_{x \in L_j}|H_i(x,y)|^{\delta s}\int_{\{(x,y) \in D_i: x \in L_j\}}|s_i(x,y)|^{s - {1 \over 2}}|(x,y)|^{-\beta}\,dx\,dy \eqno (4.31)$$
Now one may argue as after $(4.23a)$ to get the desired bound $\sum_j |N_j(\lambda)| \leq  C_s''|\lambda|^{-{1 \over 2} - \eta_s}$ for some positive 
$\eta_s$.

\noindent {\bf Case 4: $\beta_i = 0$, $\eta_i$ is not linear, and $B_{il} = 1$.}

We now suppose $\beta_i = 0$ and $B_{il} = 1$, where recall $(A_{il}, B_{il})$ denotes the upper vertex of the lowest nonhorizontal edge of the Newton polygon of $s_{il}(x,y)$, as in Lemma 2.5. Let $-{1 \over m_{il}}$ denote the slope of this lowest nonhorizontal edge. Our constructions are such that
the domain $D_{il}'$ of $(3.12b)$ must be a subset of $\{(x,y): |y| < C|x^{m_{il}}|\}$ for some $C$, roughly speaking since a term of the form
 $d_i x^{\alpha_i}$ cannot dominate the Taylor expansion of $s_{il}(x,y)$ on any larger set. The Newton polygon of $\partial_y s_{il}(x,y)$ has a 
 vertex at $(A_{il},0)$ which is the lower vertex of an edge of this Newton polygon of slope less than $-{1 \over m_{il}}$. Thus $|\partial_y s_{il}(x,y)|
 \sim x^{A_{il}}$ on $D_{il}'$. By $(2.2)$ one also has $|\partial_{xy} s_{il}(x,y)| \sim x^{A_{il} - 1}$ on $D_{il}'$.
 
 We decompose $(3.13b)$ dyadically in the $x$ variable only this time. Again let $L_j$ be the interval $[2^{-j-1},2^j]$. We will bound $\sum_j|M_j(\lambda)|$,
  where $M_j(\lambda)$ is given by
 $$\int_{\{(x,y) \in D_{il}' \cap L_j:\,|H_{il}(x,y)| < |\lambda|^{-{1 \over 100}}{\,\, {\rm or}\,\,}|\eta_{il}(x,y)| < |\lambda|^{-{1 \over 100}}\}} e^{-i\lambda_0 x - i\lambda_2y - \lambda_2K_{il}(x) - i\lambda_3s_{il}(x,y)} $$
$$\times |s_{il}(x,y)|^z |H_{il}(x,y)|^{\delta z} e^{z^2} \rho_{il}(x,y)\,dx\,dy \eqno (4.32)$$
We apply Lemma 3.3, using that $|\partial_{xy} s_{il}(x,y)| \sim 2^{-j(A_{il} - 1)}$. The $N$ of $(3.5b)$ is determined as in the $\beta_i = 1$ case, and is
given by
$$C_s \sup_{L_j}|H_{il}(x,y)|^{\delta s} \times 2^{-js\alpha_i}  \times  2^{-j\beta}$$
So if we let $a_{il}$ be such that the width of the vertical cross sections of $D_{il}'$ is 
$\sim x^{a_{il}}$, and recall that $|s_{il}(x,y)| \sim x^{\alpha_i}$ here, then Lemma 3.3 gives that
$$|M_j(\lambda)| \leq  C_s {1 \over (|\lambda| 2^{-j(A_{il} -1)})^{1 \over 2}} \sup_{x \in L_j}|H_{il}(x,y)|^{\delta s} \times 2^{-js\alpha_i} \times  2^{-j\beta} \times 
2^{-{j \over 2}} \times 2^{-{a_{il} \over 2}j} \eqno (4.33)$$
Since $D_{il}' \subset \{(x,y): |y| < C|x^{m_{il}}|\}$, we must have $a_{il} \geq m_{il}$. Furthermore, since $(\alpha_i,0)$ and $(A_{il},1)$ are joined
by an edge of the Newton polygon of $s_{il}$ of slope $-{1 \over m_{il}}$, we also have $\alpha_i - A_{il} = m_{il}$. Thus $(4.33)$ implies that
$$|M_j(\lambda)| \leq  C_s {1 \over (|\lambda| 2^{-j(\alpha_i - m_{il}  -1)})^{1 \over 2}} \sup_{x \in L_j}|H_{il}(x,y)|^{\delta s} \times 2^{-js\alpha_i} \times  2^{-j\beta} \times 
2^{-{j \over 2}} \times 2^{-{m_{il} \over 2}j}$$
$$=  C_s {1 \over (|\lambda| 2^{-j \alpha_i})^{1 \over 2}} \sup_{x \in L_j}|H_{il}(x,y)|^{\delta s} \times 2^{-js\alpha_i} \times  2^{-j\beta} \times 
2^{-j} \times 2^{-m_{il}j}\eqno (4.34)$$
There are $c, r > 0$ such that $|s_{il}(x,y)| \sim x^{\alpha_i}$ on $\{(x,y): 0 < x < r, 0 < y < c x^{m_{il}}\}$. Thus we may rewrite $(4.34)$ in the
more convenient form
$$|M_j(\lambda)| \leq  C_s' |\lambda|^{-{1 \over 2}}\sup_{x \in L_j}|H_{il}(x,y)|^{\delta s}\int_{\{(x,y): x \in L_j, \, 0 <  y < c x^{m_{il}}\}}|s_{il}(x,y)|^{s - {1 \over 2}}|(x,y)|^{-\beta}\,dx\,dy \eqno (4.35)$$
This is the estimate of $(4.31)$, so once again one may argue as after $(4.23a)$ to get the desired bound $\sum_j |M_j(\lambda)| \leq  C_s''|\lambda|^{-{1 \over 2} - \eta_s}$ for some positive $\eta_s$.

\noindent {\bf Case 5: $\beta_i = 0$, $\eta_i$ is not linear, and $B_{il} > 1$.}

This is the most difficult case. Once again we are bounding $(3.13b)$. Note that this time we use $|s_{il}^*(x,y)|^z$ in the damping function for $|\bar{s}_{il}(x,y)|^z$. The second $x$ derivative of the phase in $(3.13b)$ is given by $-\lambda_2 K_{il}''(x) - \lambda_3 
\partial_{xx} s_{il}(x,y)$. By $(2.1)$, for some $d_{il} > 0$ and a small $\eta > 0$ the term $\lambda_3 \partial_{xx} s_{il}(x,y)$ is between
$\lambda_3(1 - \eta)d_{il}\alpha_i(\alpha_i - 1) x^{\alpha_i - 2}$ and $\lambda_3(1 + \eta)d_{il} \alpha_i(\alpha_i - 1)  x^{\alpha_i - 2}$. The
 Taylor expansion of $K_{il}(x)$ is of the form 
$l_{il}x^{s_i} + o(x^{s_i})$, where $l_{il} \neq 0$ and $s_i > 1$. We write $s_i$ and not $s_{il}$ because this exponent is the same for all $l$, given $i$.

 Thus on a small enough neighborhood of $x = 0$, we have that $K_{il}''(x)$ is between $(1 - \eta)l_{il}s_i(s_i - 1)x^{s_i-2}$ and 
 $(1 + \eta) l_{il}s_i(s_i - 1) x^{s_i-2}$. Consequently, using the fact that
 $1 < s_i < \alpha_i$ as was shown in [G2], on a small enough neighborhood of the origin, if $x_0$ is such that
 $- \lambda_2 l_{il} s_i(s_i - 1)x_0^{s_i - 2} - \lambda_3 d_{il}\alpha_i(\alpha_i - 1)x_0^{\alpha_i -2} = 0$, then for some $c > 0$ and some $0 < s < 1$,
  if ${x \over x_0} > {1 \over s}$ or ${x \over x_0} < s$  then one has
 $$ |-\lambda_2 K_{il}''(x) - \lambda_3 \partial_{xx} s_{il}(x,y)| > c|\lambda_3 \partial_{xx} s_{il}(x,y)| \eqno (4.36)$$
 If the signs of $\lambda_2 l_{il}$ and $\lambda_3 d_{il}$ are the same, so that no such $x_0$ exists, then $(4.36)$ will always hold, and as a result 
 one can bound $(3.13b)$ by $C_s|\lambda|^{-{1 \over 2} - \eta}$ exactly as one proved the bounds in the case when $\beta_i = 0$ and $\eta_i$ 
 was linear. If $x_0$ does exist, one can bound the ${x \over x_0} > {1 \over s}$ or ${x \over x_0} < s$ portion of $(3.13b)$ by $C_s|\lambda|^{-{1 \over 2} - \eta}$ in the
 same fashion. Because we are now using $|s_{il}^*(x,y)|^z$ in the damping function, we do need to make use of the fact that the integral of $(s_{il}^*(x,y))^{-t}$ over $D_{il}'$ is finite whenever $t < g$, which follows from Lemma 2.6. 
 
 Thus in our future arguments it suffices to prove our bounds for the $s < {x \over x_0} < {1 \over s}$ portion only.
Before proceeding any further, we must do a further coordinate change on those $D_{il}'$ whose lower boundaries are not on the $x$-axis. Using
a version of the notation of Theorem 2.2, we write  the equation of this boundary as $y = h_{il}x^{m_{il}}
+ o(x^{m_{il}})$ for some $h_{il} > 0$. We do a coordinate change turning $(x,y)$ into  $(x, y + (h_{il} - \mu)x^{m_{il}})$ for some small $\mu$ 
determined by the following conditions. By Theorem 2.2, there are $(a_{il},b_{il})$ and $(\gamma_{il},\delta_{il})$
such that  $|\partial_y s_{il}(x,y)| \sim x^{a_{il}}y^{b_{il}}$ and $|\partial_{yy} s_{il}(x,y)| \sim x^{\gamma_{il}}y^{\delta_{il}}$ on $D_{il}$. 
The lower boundary of the transformed $D_{il}'$ is of the form $= \mu x^{m_{il}} + o(x^{m_{il}})$. Let $\tilde{s}_{il}(x,y)$ be the transformed 
$s_{il}(x,y)$ in the new coordinates. The constant $\mu$ is chosen to be small enough
so that for some $r > 0$, for all $0 < x < r$ and  $0 \leq  y < \mu  x^{m_{il}}$ we have
$$|\tilde{s}_{il}(x,y)| \sim x^{\alpha_i} {\hskip 0.5 in}|\partial_y \tilde{s}_{il}(x,y)| \sim  x^{a_{il}}x^{m_{il}b_{il}} {\hskip 0.5 in} |\partial_{yy} \tilde{s}_{il}(x,y)| \sim x^{\gamma_{il}}x^{m_{il}\delta_{il}} \eqno (4.37)$$
  That we can ensure that $(4.37)$ holds follows from the
fact that in the proofs of Theorem 2.1 and 2.2 there is some slack in the sense that $h_{il}$ can be replaced by $h_{il} - \mu$ for sufficiently
small $\mu > 0$.

Note that in the new coordinates we still have $|\tilde{s}_{il}(x,y)| \sim x^{\alpha_i}$, $|\partial_y \tilde{s}_{il}(x,y)| \sim x^{a_{il}}y^{b_{il}}$, and $|\partial_{yy} \tilde{s}_{il}(x,y)| \sim x^{\gamma_{il}}y^{\delta_{il}}$. This is because $y$ in the old coordinates is comparable in magnitude to $y$ in the new coordinates. Also, $(\alpha_i,0)$ is the lowest vertex of the Newton polygon of $\tilde{s}_{il}(x,y)$, the $B_{il}$ and $B_{il}'$ of Lemmas 2.5 and 2.6 are unchanged, and those lemmas will still hold in the new coordinates,
again using that $y$ in the old coordinates is comparable in magnitude to $y$ in the new coordinates.
Using the chain rule, we also see that Corollary 2.3 still holds. In fact, the bounds of Corollary 2.3 will also hold
when $0 \leq y < \mu x^{m_{il}}$ if $\mu$ were chosen small enough.

So for the purpose of our subsequent arguments, which will not use any aspects of Theorem 2.1 and 2.2 that do not hold in the new coordinates, we are able to replace $s_{il}(x,y)$ by $\tilde{s}_{il}(x,y)$, replace $|s_{il}^*(x,y)|^z$ by the 
corresponding $|\tilde{s}_{il}^*(x,y)|^z$, replace $D_{il}'$ by the transformed domain, and add the condition $(4.37)$. $K_{il}(x)$ 
gets replaced by $K_{il}(x) + (h_{il} - \eta)x^{m_{il}}$. In view of the defintion of $x_0$, one might ask if $x_0$ in the old coordinates can be shifted
enough in the new coordinates so that the $s < {x \over x_0} < {1 \over s}$ argument can be dealt with like the ${x \over x_0} < s$ or 
${x \over x_0} > {1 \over s}$ situation was dealt with above. 
Unfortunately that is not necessarily the case; it turns out that typically $m_{il}$ is large enough in the cases being considered that $x_0$ remains
unchanged. 

We now proceed to bound the $s < {x \over x_0} < {1 \over s}$ portion of $(3.13b)$.  
Let $e_{il}x^{\alpha_i}$ be the initial term of the Taylor expansion of $s_{il}(x,0)$. We define $x_1$ by the condition that
$$\lambda_2 l_{il} s_i(s_i - 1)x_1^{s_i - 2} + \lambda_3 e_{il} \alpha_i(\alpha_i - 1)x_1^{\alpha_i -2} = 0 \eqno (4.38)$$
Here as before $K_{il}(x) = l_{il}x^{s_i} + o(x^{s_i})$.
Since $x_0$ was defined similarly to $x_1$, with $e_{il}$ replaced by $d_{il}$, for some constant $C$ one has ${1 \over C}x_1 < x_0 < Cx_1$. Thus
the $(x,y)$ for which $s < {x \over x_0} < {1 \over s}$ is a subset of the $(x,y)$ for which $s' < {x \over x_1} < {1 \over s'}$ for some $0 < s' < s$. 
 Similarly, we define $x_2$ by the condition that
$$\lambda_2 \partial_{xx}s_{il}(x_2,0) + \lambda_3 K_{il}''(x_2) = 0 \eqno (4.39)$$
Since  $l_{il}s_i(s_i - 1)x^{s_i - 2}$ is the leading term of $\partial_{xx}s_{il}(x,0)$ and $e_{il}\alpha_i(\alpha_i - 1)x^{\alpha_i}$ is the leading term of  $K_{il}''(x)$, on a small enough neighborhood of the origin we have that the points where $s < {x \over x_0} < {1 \over s}$ is a subset of the points where
$t < {x \over x_2} < {1 \over t}$, where $0 < t < s'$. Thus for our purposes it suffices to
bound the portion of $(3.13b)$ for which $t < {x \over x_2} < {1 \over t}$. 

\begin{lemma} On a sufficiently small neighborhood of the origin, at the points where $t < {x \over x_2} < {1 \over t}$ there is a constant $C$ for
which
$$|\lambda_2 \partial_{xx}s_{il}(x,0) + \lambda_3 K_{il}''(x)| > C|\lambda||x_2^{\alpha_i-3}(x - x_2)| \eqno (4.40)$$
\end{lemma}

\noindent {\bf Proof.} First note that since $|\lambda_1| + |\lambda_2| < |\lambda_3|$ here, we can replace $|\lambda|$ by $|\lambda_3|$
when proving $(4.40)$. 
We use the fact that the first $x$ derivative of $\lambda_2 \partial_{xx}s_{il}(x,0) + \lambda_3 K_{il}''(x)$ is of the form 
 $$\lambda_2 \big( l_{il}s_i(s_i - 1)(s_i - 2)x^{s_i - 2} + o(x^{s_i - 3}) \big)+ \lambda_3 \big(e_{il}\alpha_i(\alpha_i - 1)(\alpha_i - 2)
 x^{\alpha_i- 3} + o(x^{\alpha_i - 3}) \big) \eqno (4.41) $$
 The two terms in $(4.40)$ are of comparable magnitude, but because $\alpha_i \neq s_i$ they will no longer cancel near $x = x_2$. Instead, if 
 $\delta > 0$ is sufficiently small, then on $|x - x_2| < \delta|x_2|$, $(4.41)$ will be of magnitude greater than some $C$ times the 
 magnitude of the individual terms, giving a lower bound of $C|\lambda|x_2^{\alpha_i - 3}$. Given this lower bound for the magnitude of the derivative of 
 $\lambda_2 \partial_{xx}s_{il}(x,0) + \lambda_3 K_{il}''(x)$ and the fact that this function has a zero at $x = x_2$, the estimate $(4.40)$ follows
 on $|x - x_2| < \delta|x_2|$.

On the other hand, suppose $|x - x_2| > \delta|x_2|$. On a sufficiently small neighborhood of the origin, up to 
error terms the left-hand side of $(4.40)$ is given by $|\lambda_2 l_{il} s_i(s_i - 1)x^{s_i - 2} + \lambda_3 e_{il} \alpha_i(\alpha_i - 1)x^{\alpha_i -2}|$,
which has a zero at $x = x_1$. If we are near enough to the origin, then $|x_2/x_1|$ will be between $1 - {\delta \over 4}$ and $1 + {\delta \over 4}$ and
$|x - x_2| > \delta|x_2|$ will imply $|x - x_1| > {\delta \over 2}|x_1|$, which will imply that $|\lambda_2 l_{il} s_i(s_i - 1)x^{s_i - 2} + \lambda_3 e_{il} \alpha_i(\alpha_i - 1)x^{\alpha_i -2}|$ is bounded below by an expression of the form of the right-hand side of $(4.40)$. If we are in a sufficiently 
small neighborhood of the origin, then the error terms will be small enough so that $(4.40)$ holds for $t < {x \over x_2} < {1 \over t}$. This completes the 
proof of Lemma 4.1.
 
 We now divide the integral $(3.13b)$ dyadically in $x$ and $y$, centered at $(x_2,0)$, recalling that $\bar{s}_{il}(x,y) = s_{il}^*(x,y)$ now, where
 $s_{il}^*(x,y)$ is as in Lemma 2.6. Namely,
let $L_{jk} = \{(x,y): 2^{-j -1} < |x - x_2| \leq 2^{-j},\,2^{-k-1} < y < 2^{-k}\}$, and we seek bounds for $\sum_{jk} |I_{jk}|$, where
$$I_{jk} = \int_{\{(x,y) \in D_{il}' \cap L_{jk}: \,|H_{il}(x,y)| < |\lambda|^{-{1 \over 100}}{\,\, {\rm or}\,\,}|\eta_{il}(x,y)| < |\lambda|^{-{1 \over 100}}\}} e^{-i\lambda_0 x - i\lambda_2y - \lambda_2K_{il}(x) - i\lambda_3s_{il}(x,y)} $$
$$\times |s_{il}^*(x,y)|^z |H_{il}(x,y)|^{\delta z} e^{z^2} \rho_{il}(x,y)\,dx\,dy \eqno (4.42)$$
Here $H_{il}(x,y)$, $\eta_{il}(x,y)$, and so on are replaced by the corresponding functions if we needed to do the additional coordinate change 
described earlier. The key estimate for $|I_{jk}|$ is provided by the following lemma.

\begin{lemma} For some constant $C$ one has the following, where as before $s = Re\,z$.
$$|I_{jk}| \leq C \sup_{D_{il}' \cap L_{jk}}|H_{il}(x,y)|^{\delta s} \int_{L_{jk}} {|s_{il}^*(x,y)|^s  \over (|\lambda| x^{\alpha_i - 3}|x - x_2|^3)^{1 \over 2} + (|\lambda| x^{\gamma_{il}}y^{\delta_{il} + 2})^{1 \over 2}}|(x,y)|^{-{\beta}}\,dx\,dy \eqno (4.43)$$
\end{lemma}

\noindent {\bf Proof.} The second $y$ derivative of the phase in $(4.42)$ is given by $\lambda_3 \partial_{yy}s_{il}(x,y)$, which is of magnitude
$\sim |\lambda_3|x^{\gamma_{il}}y^{\delta_{il}} \sim |\lambda|x^{\gamma_{il}}y^{\delta_{il}}$. One can apply the argument for the
$\beta_i \geq 2$ case dealt with earlier, and the result is the analogue of $(4.23a)$, which will be $(4.43)$ with a denominator of
$(|\lambda| x^{\gamma_{il}}y^{\delta_{il} + 2})^{1 \over 2}$ in place of $ (|\lambda| x^{\alpha_i - 3}|x - x_2|^3)^{1 \over 2} + $
$ (|\lambda| x^{\gamma_{il}}y^{\delta_{il} + 2})^{1 \over 2}$. Thus in order to prove $(4.43)$ one need only consider the case where for some 
small constant $\delta_0$, on $L_{jk}$ one has
$$  x^{\alpha_i - 3}|x - x_2|^3 > {1 \over \delta_0} x^{\gamma_{il}}y^{\delta_{il} + 2}  \eqno (4.44)$$
The constant $\delta_0$ will be determined by our arguments. We write the Taylor expansion of $\partial_y s_{il}(x,0)$ as 
$n_{il}x^{p_{il}} + o(x^{p_{il}})$ if the expansion is nonzero. We focus now on the case where this Taylor expansion is  not identically zero and where
 on $L_{jk}$ we have
$$ x^{\alpha_i - 3}|x - x_2|^3 <  {1 \over \delta_0^{1 \over 2}} \,x^{p_{il}-1}|x - x_2| y \eqno (4.45)$$ 
Then in view of $(4.44)$, on $L_{jk}$ we have
$$ x^{\gamma_{il}}y^{\delta_{il} + 2} < \delta_0^{1 \over 2} x^{p_{il}-1}|x - x_2|y   \eqno (4.46)$$
Next, observe that 
$$\partial_{xy}s_{il}(x,y) = \partial_{xy} s_{il}(x,0) + \int_0^y \partial_{xyy}s_{il}(x,z)\,dz \eqno (4.47)$$
Note that by Corollary 2.3, which after the
final coordinate change applies down to $y = 0$, we have
$$\bigg|\int_0^y \partial_{xyy}s_{il}(x,z)\,dz\bigg| \leq C{1 \over x}\bigg|\int_0^y \partial_{yy}s_{il}(x,z)\,dz\bigg| \eqno (4.48)$$
Using $(4.37)$ if necessary, the integral on the right is bounded $C'x^{\gamma_{ij} - 1}y^{\delta_{ij}+1}$, which by $(4.46)$ is in turn bounded by
$C\delta_0 x^{p_{il} - 2}|x - x_2|$. Since $|x - x_2| < C|x|$ due to the fact that $t < {x \over x_2} < {1 \over t}$, this is in turn bounded by $C'\delta_0 x^{p_{il} - 1}$.

 On the other hand, by the definition of $p_{il}$ one has $ \partial_{xy} s_{il}(x,0) = n_{il}p_{il}x^{p_{il}-1} + o(x^{p_{il}-1})$.  So by inserting 
 the above bounds for the integral back in $(4.47)$, if $\delta_0$ were chosen appropriately small one gets
$$|\partial_{xy}s_{il}(x,y)| >\bigg| {n_{il} p_{il} \over 2}x^{p_{il} - 1}\bigg| \eqno (4.49)$$ 
So using $(4.45)$, for some constant $C$, on
$L_{jk}$ we have
$$|\partial_{xy}s_{il}(x,y) (x - x_2) y| >  C\delta_0^{1 \over 2}x^{\alpha_i - 3}|x - x_2|^3 \eqno (4.50)$$ 
Since the $xy$ derivative of the phase function in $(4.42)$ is exactly $\lambda_3 \partial_{xy}s_{il}(x,y)$, we can now apply the argument of the $\beta_i = 1$ situation given above, which led to $(4.23a)$. This time, we get the bound
$$|I_{jk}| \leq C \sup_{D_{il}' \cap L_{jk}}|H_{il}(x,y)|^{\delta s} \int_{L_{jk}} {|s_{il}^*(x,y)|^s  \over |\lambda\partial_{xy}s_{il}(x,y) (x - x_2) y| ^{1 \over 2}}|(x,y)|^{-{\beta}}\,dx\,dy \eqno (4.51)$$
In view of $(4.50)$ and $(4.44)$, one can replace the denominator in $(4.51)$ by the denominator in $(4.43)$. Thus we are done with the proof of 
Lemma 4.1 in the case where $(4.44)$ and $(4.45)$ hold.

It remains to consider the possibility that $(4.44)$ holds but $(4.45)$ does not, or that the Taylor expansion of $\partial_y s_{il}(x,0)$ is identically zero.
In the former case, since $(4.45)$ does not hold everywhere on $L_{jk}$, there is a constant
 $C$ such that on $L_{jk}$ one has
$$x^{p_{il}-1}|x - x_2| y < C \delta_0^{1 \over 2}\,x^{\alpha_i - 3}|x - x_2|^3 \eqno (4.52)$$
We will apply the Van der Corput lemma for second derivatives in the $x$ direction. We take the second $x$ derivative of the phase function
in $(4.42)$ and Taylor expand the $\lambda_3\partial_{xx}s_{il}(x,y)$ term in $y$, resulting in 
$$-\lambda_3 \partial_{xx} s_{il}(x,0) -\lambda_2K_{il}''(x) -\lambda_3\partial_{xxy}s_{il}(x,0)y +  \lambda_3\int_0^yz\partial_{xxyy}s_{il}(x,z)\,dz
\eqno (4.53)$$
The idea now is that the sum of the first two terms in $(4.53)$  are of absolute value at least $C|\lambda||x_2^{\alpha_i - 3}(x- x_2)|$ by Lemma 
4.1, and the magnitude of the remaining terms will be much smaller by $(4.44)$ and $(4.52)$. We first look at the $-\lambda_3\partial_{xxy}s_{il}(x,0)y$
term. By Corollary 2.3, which in the final coordinates will apply on the $x$ axis, if this term is nonzero we have that 
$$ |\partial_{xxy}s_{il}(x,0)| \leq C{1 \over x^2} |\partial_{y}s_{il}(x,0)| \leq C' x^{p_{il}-2} \eqno (4.54)$$
Thus in view of $(4.52)$ we have
$$|\lambda_3\partial_{xxy}s_{il}(x,0)y| \leq C''\delta_0^{1 \over 2}|\lambda||x^{\alpha_i - 4}(x - x_2)^2 \eqno (4.55)$$
Given that $|x - x_2| < C|x|$ whenever $t < {x \over x_2} < {1 \over t}$, this implies that
$$|\lambda_3\partial_{xxy}s_{il}(x,0)y| \leq C'''\delta_0^{1 \over 2}|\lambda||x^{\alpha_i - 3}|x - x_2| \eqno (4.56)$$
Thus if $\delta_0$ is sufficiently small, $|\lambda_3\partial_{xxy}s_{il}(x,0)y |$ will be small in comparison to the magnitude of the sum of the first
two terms in $(4.53)$. 

\noindent Proceeding now to the integral term of $(4.53)$, by Corollary 2.3 we have
$$\bigg|\int_0^y z\partial_{xxyy}s_{il}(x,z)\,dz\bigg| \leq C{1 \over x^2} \int_0^y z|\partial_{yy}s_{il}(x,z)|\,dz \eqno (4.57)$$
Using that $|\partial_{yy}s_{il}(x,y)| \sim| x^{\gamma_{il}}y^{\delta_{il}}|$ on $D_{il}'$ and using $(4.37)$ below $D_{il}'$ if necessary, we have
$$C{1 \over x^2} \int_0^y z|\partial_{yy}s_{il}(x,z)|\,dz \leq C'x^{\gamma_{il}- 2}y^{\delta_{il}+2} \eqno (4.58)$$
Using $(4.44)$, this is bounded by $C'\delta_0\, x^{\alpha_i - 5}|x - x_2|^3$. Using again that $|x - x_2| < C|x|$, this is in turn bounded by
$C''\delta_0\, x^{\alpha_i - 3}|x - x_2|$. Thus if $\delta_0$ were chosen sufficiently small, the integral term in $(4.53)$ is also of much smaller magnitude
than the sum of the first two terms. 

Putting the above together, we see that if $\delta_0$ were chosen appropriately small, the sum of the first two terms dominate $(4.53)$, and thus this second derivative of the phase in $(4.42)$ is
bounded below by $Cx^{\alpha_i - 3}|x - x_2|$. We now apply Lemma 3.2 for second derivatives in the $x$-direction, similar to how we did in the
$\beta_i = 0$, $\eta_i(x,y)$ linear case. The result is the bound
$$|I_{jk}| \leq C \sup_{D_{il}' \cap L_{jk}}|H_{il}(x,y)|^{\delta s} \int_{L_{jk}} {|s_{il}^*(x,y)|^s  \over |\lambda x^{\alpha_i - 3}(x - x_2)^3|^{1 \over 2}}|(x,y)|^{-{\beta}}\,dx\,dy \eqno (4.59)$$
Since $(4.44)$ is assumed to hold, $(4.43)$ is therefore satisfied. This completes the proof of Lemma 4.2 for the situation where the Taylor expansion of $\partial_ys_{il}(x,0)$ is not identically zero, and $(4.44)$ holds but $(4.45)$ does not.

It remains to consider the case where $(4.44)$ holds and the Taylor expansion of $\partial_ys_{il}(x,0)$ is identically zero. In this case, in place of 
$(4.53)$ we have that the second $x$ derivative of the phase function is now given by
$$-\lambda_3 \partial_{xx} s_{il}(x,0) -\lambda_2K_{il}''(x) +  \lambda_3\int_0^yz\partial_{xxyy}s_{il}(x,z)\,dz
\eqno (4.60) $$
Then exactly as after $(4.58)$, the integral in $(4.60)$ is of far smaller magnitude than the magnitude of the sum of the first two terms if $\delta_0$ were
chosen appropriately small. So once again the second $x$ derivative of the phase is bounded below by $Cx^{\alpha_i - 3}|x - x_2|$ and $(4.59)$ holds
exactly like before. This concludes the proof of Lemma 4.2.

By Lemma 2.5, we have $x^{\gamma_{il}}y^{\delta_{il} + 2} > Cy^{B_{il}'}$ where $B_{il}'$ is in that lemma. As a result, using the arithmetic geometric
mean inequality, the denominator of $(4.43)$ satisfies
$$(|\lambda| x^{\alpha_i - 3}|x - x_2|^3)^{1 \over 2} + (|\lambda| x^{\gamma_{il}}y^{\delta_{il} + 2})^{1 \over 2} > 
C(|\lambda| x^{\alpha_i - 3}|x - x_2|^3)^{1 \over 2} + (|\lambda| y^{B_{il}'})^{1 \over 2}$$
$$\geq 2C|\lambda|^{1 \over 2} x^{\alpha_i - 3 \over 4}|x - x_2|^{3 \over 4} y^{B_{il}' \over 4} \eqno (4.61)$$
Note that the right-hand side of $(4.61)$ is exactly $2C|\lambda|^{1 \over 2}(s_{il}^*(x,y))^{1 \over 2}x^{-{3 \over 4}}|x - x_2|^{3 \over 4}$. Thus $(4.43)$ implies
that
$$|I_{jk}| \leq C |\lambda|^{-{1 \over 2}}\sup_{D_{il}' \cap L_{jk}}|H_{il}(x,y)|^{\delta s} \int_{L_{jk}} {(s_{il}^*(x,y))^s  \over (s_{il}^*(x,y))^{1 \over 2}x^{-{3 \over 4}}|x - x_2|^{3 \over 4}}|(x,y)|^{-{\beta}}\,dx\,dy \eqno (4.62)$$
We now add $(4.62)$ over all $j$ for fixed $k$. Letting $M_k = \{(x,y): t < {x \over x_2} < {1 \over t},\, 2^{-k-1} < y \leq 2^{-k}\}$, we get 
$$\sum_j |I_{jk}| \leq C|\lambda|^{-{1 \over 2}} \sup_{D_{il}' \cap M_k}|H_{il}(x,y)|^{\delta s}   \int_{M_k}{(s_{il}^*(x,y))^s  \over  (s_{il}^*(x,y))^{1 \over 2}x^{-{3 \over 4}}|x - x_2|^{3 \over 4}}|(x,y)|^{-{\beta}}\,dx\,dy  \eqno (4.63)$$
Since $x$ appears to the $-{3 \over 4}$ power in $(4.63)$, examining the $x$ integral in $(4.63)$ for fixed $y$ we see that
$$\sum_j |I_{jk}| \leq C |\lambda|^{-{1 \over 2}}\sup_{D_{il}' \cap M_k}|H_{il}(x,y)|^{\delta s}  \int_{M_k}(s_{il}^*(x,y))^{s- {1 \over 2}} |(x,y)|^{-{\beta}}\,dx\,dy  \eqno (4.64)$$
Adding this over all $k$ gives the right-hand side of $(4.31)$. Since by Lemma 2.6 the integral of $(s_{il}^*(x,y))^{-t}$ over $D_{il}'$ is finite
 whenever $t < g$, one may argue as after $(4.23a)$ to get the desired bound of $C_s|\lambda|^{-{1 \over 2} - \nu_s}$ for $\sum_{j,k}|I_{jk}|$. This
completes the argument for the Fourier transform decay estimates when $\beta_i = 0$, $\eta_i$ is not linear, and $B_{il} > 1$.

\subsubsection {The end of the proof of Theorem 1.1 a).}

Let $d\sigma_z$ denote the damped surface measure being dilated in $(3.8)$. Adding up the estimates in the various cases above gives us an overall
bound of $|\widehat{d\sigma_z}| \leq C_s|\lambda|^{-{1 \over 2} - \nu_s}$ when $s > \max(0,{1 \over 2} - g)$, where as always $s$ denotes $Re\,z$. 
Thus by Theorem 3.1 we have $||M_z f||_{L^2} \leq C_s||f||_{L^2}$ for such $z$. In the beginning of this section, we also showed that
for $z > -g$ one has $||M_z f||_{L^{\infty}} \leq C_s||f||_{L^{\infty}}$. So in particular this holds for $s > \max(-{1 \over 2}, -g)$. 

Note that
$0 = \alpha(\max(-{1 \over 2}, -g)) + (1 - \alpha)\max(0,{1 \over 2} - g)$, where $\alpha = \max(1 - 2g, 0)$. Thus  if $1 > \beta > \alpha$,
 one can write $0 = \beta s + (1 - \beta) s'$, where $0 > s > \max(-{1 \over 2}, -g)$ and  $s' >  \max(0,{1 \over 2} - g)$. Hence by complex interpolation $M = M_0$ is bounded on $L^p$, where
 ${1 \over p} = \beta {1 \over \infty} + (1 - \beta) {1 \over 2}$. In other words $p = {2 \over 1 - \beta}$. As $\beta$ approaches $\alpha$, this 
 exponent approaches ${2 \over 1 - \alpha} = {2 \over 1 - \max(1 - 2g,0)} = {2 \over \min(2g,1)} = \max(2, {1 \over g})$. Thus $M$ is bounded on 
$L^p$ for $p > \max(2,{1 \over g})$ as needed. This completes the proof of Theorem 1.1a).
 
 \subsubsection {The proof of Theorem 1.1 b).}
 
 We assume we are in the setting of Theorem 1.1b); that is, we assume that he tangent plane to $S$ at $(x_0,y_0,z_0)$ does not contain 
 the origin and on some a neighborhood of $(x_0,y_0)$ one has $|\phi(x,y)| > C|{\bf x} - {\bf x_0}|^{-\beta}$ for some $C > 0$. We can assume
 $g < 1$ since the $g = 1$ case occurs only when $\beta = 0$ and the surface is nondegenerate, where it is known $M$ is bounded on $L^p$ 
 only if $p > {3 \over 2}$. 
 
 Like in the rest of
 this section, we assume that we have rotated so that $(0,0,1)$ is normal to $S$ at $(x_0,y_0,z_0)$. The tangent plane condition then implies that
 $z_0 \neq 0$. We let $f(x,y,z) = \alpha(x,y,z)|z|^{-g}|\ln |z||^{-g^*}$, where $1 \geq  g^* > g$ and where $\alpha(x,y,z)$ is a 
 nonnegative compactly supported function identically equal to 1 on a disk centered at the origin. Then $f(x,y,z) \in L^{1 \over g}(\R^3)$.
 We claim that for small $r$ one has
 $$ \int_{|(x,y)| < r} |s(x,y)|^{-g}|\ln |s(x,y)||^{-g^*}|(x,y)|^{-\beta}\,dx\,dy = \infty \eqno (4.65)$$
 In other words, integrating $|z|^{-g}|\ln |z||^{-g^*}$ over the surface near the origin with respect to the measure $|(x,y)|^{\beta}\,dx\,dy$
 results in infinity. To see why $(4.65)$ holds, observe that the definition of $g$ in terms of asymptotics implies that the portion of $(4.65)$ with 
 $2^{-j - 1} < |s(x,y)| < 2^{-j}$ is at least of order $2^{gj} j^{-g^*}\times 2^{-gj} = j^{-g^*}$. Since $g^* \leq 1$, adding this over all $j$ results in
 infinity as needed.
 
  If $(x,y,z)$ is
 such that $z \neq 0$ and $sgn(z) \neq sgn(z_0)$, then there is some $t > 0$ such that  $(x,y,z) + tS$ is tangent to the $xy$ plane at the point $(x,y,z) + t(x_0,y_0,z_0)$. So 
 by $(4.65)$ we have that $Mf(x,y,z) = \infty$ for all such $(x,y,z)$ in a neighborhood of the origin.
  We conclude that $M$ is not bounded on $L^g$, completing the proof of Theorem 1.1b). 

\subsection  {The proof of Theorem 1.2.}

The argument here is a simplified version of the proof for Theorem 1.1 so we will be brief. We saw in the proof of Theorem 1.1 that $|\bar{s}(x,y)|^{-t}$ 
is integrable on a neighborhood of the origin for all $t > -g$. Since $R_z f$ is the convolution of $f$ with a surface with density $|\bar{s}(x,y)|^z$ times a factor that is
uniformly bounded for fixed $s = Re\,z$, by Young's inequality we have that for any $1 < p < \infty$ and any $s > -g$ 
we have $||R_z f||_p \leq C_{s,p}||f||_p$. In particular this holds if $s > \max (-{1 \over 2}, g)$, which is what we will use.

Comparing $(3.17)$ with $(3.9')$ we see that the only difference up to a magnitude one factor is the absence of the $|H(x,y)|^{-\delta z}$ factor in $(3.17)$. As a result, we can
argue as in section 4, with the following modifications. We do not need to consider the case where $|\lambda_1| + |\lambda_2| \leq |\lambda_3|$, $|H(x,y)| > |\lambda|^{-{1 \over 100}}$, $|(x,y)|  > |\lambda|^{-{1 \over 100}}$, and in all subsequent parts of the argument we do not 
stipulate the condition that  $|H(x,y)| > |\lambda|^{-{1 \over 100}}$ or $|(x,y)|  > |\lambda|^{-{1 \over 100}}$, only that
$|\lambda_1| + |\lambda_2| \leq |\lambda_3|$. The effect of this is that
we will end out not having the $\sup |H_i(x,y)|^{\delta s}$ factor in our various estimates, and the result of this is that instead of ending out with 
bounds of the form $C_s|\lambda|^{-{1 \over 2} - \nu_s}$ we just have bounds of the form $C_s|\lambda|^{-{1 \over 2}}$. So the end result will be
that if $s = Re\,z > \max({1 \over 2} -g,0)$, then the expression in $(3.17)$ is bounded by $C_s|\lambda|^{-{1 \over 2}}$.

Since $(3.17)$ is just the Fourier transform of the surface measure in $(3.16)$ and $R_z f$ is the convolution with this surface measure, we conclude 
that if $s > \max({1 \over 2} - g,0)$ one has  $||R_z f||_{L^2_{1 \over 2}} \leq C_s'||f||_{L^2}$. We now interpolate this analogously to 
the interpolation at the end of the proof of Theorem 1.1a), once again using that $0 = \alpha(\max(-{1 \over 2}, -g)) + (1 - \alpha)\max(0,{1 \over 2} - g)$ where $\alpha = \max(1 - 2g, 0)$. So for $1 > \beta > \alpha$ and $p^* < p$, we have $||R f||_{L^{p^*}_{1 \over p}} \leq C_{p,p^*}||f||_{L^{p^*}}$
for ${1 \over p} = \beta{1 \over \infty} + (1 - \beta){1 \over 2}$. As we saw before, as $\beta$ approaches $\alpha$ this $p$ approaches 
$\max(2,{1 \over g})$. 
So for $p > p^* >  \max(2,{1 \over g})$ one has the estimate $||R f||_{L^{p^*}_{1 \over p}} \leq  C_{p,p^*}||f||_{L^{p^*}}$.

But given the translation-invariant nature of the Radon transform operator, one can use duality and say that if ${1 \over p^*} + {1 \over (p^*)'} = 1$, then we 
also have an estimate $||R f||_{L^{(p^*)'}_{1 \over p}} \leq  C_{p,p^*}||f||_{L^{(p^*)'}}$. Since $R$ is bounded on $L^1$ and $L^{\infty}$, we can interpolate
and conclude that $||R f||_{L^q_r} \leq C_{q,r}||f||_{L^q}$ whenever $({1 \over q}, r)$ is in the closed trapezoid connecting $(0,0), ({1 \over p^*}, {1 \over p}), 
({1 \over (p^*)'}, {1 \over p})$, and $(1,0)$. Taking the union of these as $(p,p^*)$ approaches $(\max(2,{1 \over g}), \max(2,{1 \over g}))$, we get that $||R f||_{L^q_r} \leq C_{q,r}||f||_{L^q}$ whenever $({1 \over q}, r)$ is in the $y <  \min({1 \over 2}, g)$ portion of the open trapezoid or triangle with edges given by the $x$ axis, the line $y = x$, the line $y = \min({1 \over 2}, g)$, and
the line $y = 1 - x$. Since the lines $y = x$ and $y = 1 - x$ join at $y = {1 \over 2}$, we can restate the line $y = \min({1 \over 2}, g)$ as simply
$y = g$. This concludes the proof of Theorem 1.2a).

As for the sharpness statement of part b), suppose we are in the setting of part b). In other words suppose that on some neighborhood of
$(x_0,y_0,z_0)$ one has $|\phi(x,y)| > C|{\bf x} - {\bf x_0}|^{-\beta}$ for some $C > 0$. Suppose $R$ were bounded from
from $L^p$ to $L^p_{\alpha}$ for some $1 < p < \infty$ and $\alpha > g$. Then by duality, $R$ would also be bounded from $L^{p'}$ to $L^{p'}_{\alpha}$
where ${1 \over p} + {1 \over p'} = 1$, so by interpolation $R$ is bounded from $L^2$ to $L^2_{\alpha}$. As a result, the Fourier transform 
of the surface measure on $S$ would decay at a rate of $C|\lambda|^{-\alpha}$ where $\alpha > g$. However, by Theorem 1.3c) of [G6]  
one can never get such a decay rate even in the normal direction. This concludes the proof of Theorem 1.2.

\section  {Generalizations to smooth surfaces.} 

We now describe how the statement and proof of Theorem 1.1 generalizes to
the case of smooth surfaces. So now we assume that $s(x,y)$ is a smooth function with $s(0,0) = 0$ and $\nabla s(0,0) = (0,0)$. If $s(x,y)$ has a zero of infinite order at the origin, our upcoming analogue of the Hessian condition of Theorem 1.1 will not be satisfied and no analogue of Theorem 1.2 will hold.
 So we assume that $s(x,y)$ has a zero of some finite order at $(0,0)$, in other words, that $s(x,y)$ is of finite type at $(0,0)$.

We first need to reformulate the definition of the index $g$, since the definition in the real analytic case relied on the asymptotics $(1.6)$, which do not
exist in the general smooth case. Instead, we base the definition on $(1.7)$. We let $g$ be the supremum of all numbers $t$ such that for any sufficiently
small $r > 0$, for all $0 < \epsilon < {1 \over 2}$ one has
$$ \int_{\{(x,y):\,|(x,y)| < r,\,|s(x,y)| < \epsilon\}} |(x,y)|^{-\beta} \,dx\,dy \leq C_{\beta, r}\epsilon^t  \eqno (5.1)$$
The fact that $s(x,y)$ is of finite type at the origin ensures that $g$ is nonzero. Note that $g$ is the supremum of all $t$ for which 
$|s(x,y)|^{-t}|(x,y)|^{-\beta}$ is integrable on a neighborhood of the origin, just as in the real analytic case.
 This is the key property of $g$ used in the earlier arguments.
With this definition of $g$ the new formulation of Theorem 1.1 is as follows.

\begin{theorem}
\

\noindent {\bf a)} Suppose the Hessian determinant of $s(x,y)$ does not vanish to infinite order at the origin.
 Then there is a neighborhood $N$ of $(x_0,y_0)$ such that if 
$\phi(x,y)$ is supported in $N$ and $(1.3)-(1.4)$ is satisfied, then $M$ is bounded on $L^p$ for $p > \max({1 \over g}, 2)$.

\noindent {\bf b)} If the tangent plane to $S$ at $(x_0,y_0,z_0)$ does not contain the origin and if on some neighborhood of $(x_0,y_0)$ one
has $|\phi(x,y)| > C|{\bf x} - {\bf x_0}|^{-\beta}$ for some $C > 0$, then $M$ is not bounded on $L^p$ for any 
$1 \leq p < {1 \over g}$.
\end{theorem}

So there are two differences in Theorem 1.1 in the smooth analogue; there is the adjustment in the Hessian condition in part a), and one no longer 
includes $p = {1 \over g}$ in the sharpness statement.

\noindent {\bf Detailed sketch of proof.}

We first have to reformulate the resolution of singularities theorems of section 2 so that they apply to smooth functions. The resolution of singularities 
theorems in [G5] apply to smooth functions, and the real analytic theorems used in this paper all derive in the end from those of [G5]. As a result, one can
prove smooth analogues to Theorems 2.1 and 2.2. The statements differ only when the lower boundary of $D_i'$ or $D_{il}'$ respectively is on the 
$x$-axis and the associated $s \circ \eta_i$ or $S_j \circ \eta_{il}$ satisfies $\beta_i > 0$. (i.e. when the monomialized function has a power of
$y$ in it). When both of these conditions occur, in place of the statements of Theorems 2.1 and 2.2, one can write $s \circ \eta_i$ or $S_j \circ \eta_{il}$ as the sum
of two functions. One function satisfies the estimates $(2.1)$ or  $(2.2)$ just as before, as well as Corollary 2.3, and the second function has a zero of
infinite order at the origin.
With this version of Theorems 2.1 and 2.2, the rest of section 2 holds as before. Lemma 2.6 holds for example since it applies only when $\beta_i = 0$
for $s \circ \eta_i$ in the application of Theorem 2.1.
 
Next, we consider the damping function that should be used in defining $M_z$ in the smooth case, as analogues to those of $(3.8)$ and $(3.16)$.
It turns out that
we have to adjust the definition of $\bar{s}(x,y)$ on those domains $D_i$ for which
$\beta_i > 2$. In the real analytic case, $\bar{s}(x,y)$ is just equal to $s(x,y)$ on these domains. What is needed for the arguments to work is a 
replacement for $\bar{s}(x,y)$  that in the coordinates of $D_i'$ grows for fixed $x$ at the same rate in $y$ as the function $x^{\alpha_i}y^{\beta_i}$, 
and which also is a function of $|\partial_{yy} s_i(x,y)|$ for fixed $x$. Because of the presence of the
additional function with a zero of infinite order at the origin in the above resolution of singularities theorems, $|s_i(x,y)|$ itself no longer works
 as $|\bar{s}_i(x,y)|$
 when $\beta_i > 2$. (It doesn't cause problems when $\beta_i = 2$.) Instead,
 in the coordinates of $D_i'$ we let $|\bar{s}(x,y)| =  x^{-2{\alpha_i \over \beta_i - 2}}|\partial_{yy} s_i(x,y)|^{\beta_i \over \beta_i - 2}$.
 
 With the new damping function, which we again denote by $\bar{s}(x,y)$, the arguments of Section 4 now can proceed with the following adjustments.
  The argument showing $L^{\infty}$ to $L^{\infty}$ boundedness of $M_z$
 when $Re\,z > g$ proceeds much as before. The key fact is that $|\bar{s}(x,y)|^{-t}|(x,y)|^{-\beta}$ is integrable on a neighborhood of the origin if $t < g$. We need
 only verify this where $\bar{s}(x,y)$ differs from before, namely on the $D_i$ where $\beta_i > 2$. But for fixed $x$, the distribution function of 
 $x^{-2{\alpha_i \over \beta_i - 2}}|\partial_{yy} s_i(x,y)|^{\beta_i \over \beta_i - 2}$ as a function of $y$ has the same growth rate that
 $x^{\alpha_i}y^{\beta_i}$ has. This can be seen for example by using the measure version of the Van der Corput lemma (see [C1]) for $\beta_i - 2$th derivatives
 on the function $\partial_{yy} s_i(x,y)$. So since $(x^{\alpha_i}y^{\beta_i})^{-t}$ is integrable on $D_i'$ for $t < g$, the same will be true of 
 $x^{-2{\alpha_i \over \beta_i - 2}}|\partial_{yy} s_i(x,y)|^{\beta_i \over \beta_i - 2}$. Thus the new $|\bar{s}(x,y)|^{-t}|(x,y)|^{-\beta}$ is indeed integrable on a neighborhood of the origin if $t < g$ and the $L^{\infty}$ to $L^{\infty}$ argument proceeds like before.

The arguments of section 4 when $|\lambda_1| + |\lambda_2| > |\lambda_3|$ and when $|\lambda_1| + |\lambda_2| \leq |\lambda_3|$, $|H(x,y)| > |\lambda|^{-{1 \over 100}}$, and $|(x,y)|  > |\lambda|^{-{1 \over 100}}$ still hold in the smooth case, so we need not concern ourselves with these 
situations. So we focus our attention to the situations where $|\lambda_1| + |\lambda_2| < |\lambda_3|$ and  $|H(x,y)| < |\lambda|^{-{1 \over 100}}$
or $|(x,y)|  < |\lambda|^{-{1 \over 100}}$. We first consider the case where  $\beta_i \geq 2$. The case where $\beta_i = 2$ proceeds as before, so we assume $\beta_i > 2$.
We adjust the arguments of that section as follows. Instead of doing a dyadic decomposition in $x$ and $y$ and add the results over the various 
rectangles, we fix $x$ and divide the $y$ intervals of integration in $(3.13a)$ into portions $J_{x,j}$ where $2^{-j-1} < |\partial_{yy} s_i(x,y)| \leq 2^{-j}$
on $J_{x,j}$. Using Lemma 3.2 for second derivatives in the $y$ direction, one shows that the integral over $J_{x,j}$ satisfies the same estimate 
that held in the real analytic case over the interval where $2^{-j-1} < |x^{\alpha_i}y^{\beta_i - 2}| \leq 2^{-j}$. One then adds this over all $j$ and
integrates the result in $x$ to achieve the same estimate as in the real analytic case.

The $\beta_i = 1$ case and the case where $\beta_i = 0$, $\eta_i$ is not linear, and $B_{il} = 1$  carry over to the smooth case, so we will not concern ourselves with these two cases here.
We next consider together the case when $\beta_i = 0$ and $\eta_i$ is linear, and the case where $\beta_i = 0$, $\eta_i$ is not linear, $B_{il} > 1$, and 
${x \over x_0} < s$ or ${x \over x_0} > {1 \over s}$ considered above $(4.36)$. These two situations can be dealt with  largely as in the real analytic situation, with one technical
issue arising. The functions $H_i(x,y)$ or $H_{il}(x,y)$, as well as their first $x$ derivatives, may no longer have
a number of zeroes in the $x$ variable that is bounded in $y$. This can cause issues when applying Lemma 3.2 when integrating back
$x$ derivatives of the damping function when bounding the integral on the right-hand side of $(3.4)$. Similar issues arise due to the condition 
$|H(x,y)| < |\lambda|^{-{1 \over 100}}$ in the domain of integration.

This issue can be solved by dividing up the 
domain of integration into squares of diameter $|\lambda|^{-\delta_0}$ for sufficiently small $\delta_0 > 0$, and then for some large $N$ 
replacing the $H_i(x,y)$ or 
$H_{il}(x,y)$ in the damping function by the sum of the first $N$ terms of its Taylor expansion at the square's center. One makes the same replacement in
the condition $|H(x,y)| < |\lambda|^{-{1 \over 100}}$ in the domain of integration.
If $N$ is large enough, given
$\delta_0$, then the difference between the original and adjusted integrals can be made less than say $C|\lambda|^{-1}$. Then in the adjusted 
integral on a given square, one performs the arguments of the real-analytic case and adds over all squares. If $\delta_0$ is small enough, the 
additional $C|\lambda|^{2\delta_0}$ one incurs from the addition will not be enough to erase the $\nu_s$ in the estimate in the bound of 
$C_s|\lambda|^{-{1 \over 2} - \nu_s}$ one obtains in the real-analytic case. This was the strategy the author took in the earlier paper [G7].

It remains to examine the situation where $\beta_i = 0$, $\eta_i$ is not linear, $B_{il} > 1$, and $|x - x_0| < s|x_0|$. The arguments of
the real analytic case do not immediately carry over because on $D_{il}'$, the function $\partial_y s_{il}(x,y)$ is  the sum of a main term
comparable to $x^{a_{il}}y^{b_{il}}$ and an error term with a zero of infinite order at the origin, and similarly 
$\partial_{yy} s_{il}(x,y)$ is  the sum of 
a main term comparable to $x^{\gamma_{il}}y^{\delta_{il}}$ and such an error term. Since $b_{il}$ and $\delta_{il}$ may be nonzero, the error
terms may interfere with the arguments from before. 

These issues only occur on those $D_{il}'$ whose lower boundary is on the $x$ axis since otherwise
the error terms for $\partial_y s_{il}(x,y)$ and $\partial_{yy}s_{il}(x,y)$ are much smaller than their corresponding main terms on $D_{il}'$
. For the same reason, if  $D_{il}'$ has its lower boundary is on the $x$
axis, then given $N$ and $c$ there exists an $r$ such that the earlier arguments will work for the portion of $D_{il}'$ outside the sliver
$\{(x,y): 0 < x < r, 0 < y < cx^N\}$. So we may restrict our consideration to the portion of $D_{il}'$ inside such a sliver.
 
But on the $x$-axis in the coordinates of such a $D_{il}'$, one has $s_{il}(x,y) \sim x^{\alpha_i}$.  Let $C_{il}$ denote the curve in the very original
coordinates that corresponds to the portion of the closure of $D_{il}'$ on the $x$-axis. So $s(x,y) \sim  x^{\alpha_i}$ on $C_{il}$. Because $\alpha_i \geq 2$, there is a small rotation
$R$, a $p\leq \alpha_i - 1$, and a $q \leq \alpha_i - 2$  for which $\partial_y (s \circ R^{-1}) \sim x^p$ and 
$\partial_{yy}(s \circ R^{-1}) \sim x^q$ on  the curve $R(C_{il})$ on a small enough neighborhood of the origin. Furthermore, if 
$c$ is small enough and $N$ is large enough, we will also have  $\partial_y (s \circ R^{-1}) \sim x^p$ and 
$\partial_{yy}(s \circ R^{-1}) \sim x^q$ on the points within vertical distance $cx^N$ of $R(C_{il})$. Thus if we perform the arguments of
this paper for $s \circ R^{-1}$  on just this sliver within vertical distance $cx^N$ of $R(C_{il})$, we will never be in a situation where the error
terms for $\partial_y (s \circ R^{-1})$ and $\partial_{yy}(s \circ R^{-1})$ cause any issues; since $\partial_y (s \circ R^{-1}) \sim x^p$ and $\partial_{yy}(s \circ R^{-1}) \sim x^{q}$ the analogue of $(a_{il},b_{il})$ is just $(p,0)$ and the analogue of $(\gamma_{il},\delta_{il})$
is just $(q,0)$. In particular, the situations $b_{il} \neq  0$ or $\delta_{il} \neq 0$ do not occur. This concludes our detailed sketch of the proof of part a) of
Theorem 5.1.

As for the sharpness statement of part b), we essentially use the same example as in the proof of Theorem 1.1b). If $p < {1 \over g}$, we let 
$h = {1 \over p}$ and then use $f(x,y,z) = \alpha(x,y,z)|z|^{-h}|\ln |z||^{-2h}$, where $\alpha(x,y,z)$ is again some nonnegative compactly
supported function identically equal to 1 on a  neighborhood of the origin. Then $f(x,y,z) \in L^p(\R^3)$, and since $h > g$ the definition of $g$ implies that
$$ \int_{|(x,y)| < r} |s(x,y)|^{-h}|\ln |s(x,y)||^{-2h}|(x,y)|^{-\beta}\,dx\,dy = \infty \eqno (5.2)$$
This can again be seen by adding over $j$ the portion of the integral over sets where $|s(x,y)| \sim 2^{-j}$. Like before, this implies $Mf = \infty$
on a set of positive measure and therefore $M$ is not bounded on $L^p$.

As for Theorem 1.2, the statement remains unchanged in the smooth case, and the modifications in the proof for the real-analytic case are largely as above.
 The only difference is that the Hessian determinant does
not appear in the proof, so the technical modifications above regarding the Hessian determinant do not have to be  made.
  
\section {References.}

\noindent [AGuV] V. Arnold, S Gusein-Zade, A Varchenko, {\it Singularities of differentiable maps
Volume II}, Birkhauser, Basel, 1988.  \setlength{\parskip}{0.3 em}

\noindent [B] J. Bourgain, {\it Averages in the plane over convex curves and maximal operators},
J. Anal. Math. {\bf 47} (1986), 69--85.

\noindent [BrNW] J. Bruna, A. Nagel, and S. Wainger,  {\it Convex hypersurfaces and Fourier transforms}, Ann. of
Math. (2), {\bf 127} (1988), no. 2, 333-365.

\noindent [C1] M. Christ, {\it Hilbert transforms along curves. I. Nilpotent groups}, Annals of Mathematics (2) {\bf 122} (1985), no.3, 575-596.

\noindent [C2] M. Christ, {\it Failure of an endpoint estimate for integrals along curves} in Fourier analysis and partial differential equations (Miraflores de la Sierra, 1992), 163-168, Stud. Adv. Math., CRC, Boca Raton, FL, 1995. 

\noindent [CNSW] M. Christ, A. Nagel, E. M. Stein, and S. Wainger, {\it Singular and maximal Radon transforms: analysis and geometry},
Ann. of Math. (2) {\bf 150} (1999), no. 2, 489-577. 

\noindent [CoMa] M. Cowling and G. Mauceri, {\it Inequalities for some maximal functions. II}, Trans. Amer. Math.
Soc., {\bf 296} (1986), no. 1, 341-365.

\noindent [dBvdE] M. de Bondt, A. van den Essen, {\it Singular Hessians}, Journal of Algebra {\bf 282} (2004) 195-204.

\noindent [D] J.J. Duistermaat, {\it Oscillatory integrals, Lagrange immersions, and unfolding of \hfill \break singularities}, Comm. Pure Appl.
Math., {\bf 27} (1974), 207-281.

\noindent [Gra] L. Grafakos, {\it Endpoint bounds for an analytic family of Hilbert transforms}, Duke Math. J. {\bf 62} (1991), no. 1, 23-59. 

\noindent [G1] M. Greenblatt, {\it Uniform bounds for Fourier transforms of surface measures in $\R^3$ with nonsmooth density},
Trans. Amer. Math. Soc. {\bf 368} (2016), no. 9, 6601-6625.

\noindent [G2] M. Greenblatt, {\it Estimates for Fourier transforms of surface measures in $\R^3$ and PDE applications}, Revista Math. 
Ibero. {\bf 32} (2016), no. 2, 419-446.

\noindent [G3] M. Greenblatt, {\it Convolution kernels of 2D Fourier multipliers based on real analytic functions}, to appear, J. Geom. Anal.

\noindent [G4] M. Greenblatt, {\it Maximal averages over hypersurfaces and the Newton polyhedron}, J. Funct. Anal. {\bf 262} (2012), no. 5, 2314-2348. 

\noindent [G5] M. Greenblatt, {\it Resolution of singularities in two dimensions and the stability of integrals}, Adv. Math., 
{\bf 226} no. 2 (2011), 1772-1802.

\noindent [G6] M. Greenblatt, {\it Fourier transforms of irregular local hypersurface measures in $\R^3$}, submitted.

\noindent [G7] M. Greenblatt, {\it $L^p$ boundedness of maximal averages over hypersurfaces in $\R^3$}, Trans. Amer. Math. Soc. {\bf 365}
 (2013), no. 4, 1875-1900. 

\noindent [Gr] A. Greenleaf, {\it Principal curvature and harmonic analysis}, Indiana Univ. Math. J. {\bf 30} (1981), no. 4, 519--537.

\noindent [IKeM] I. Ikromov, M. Kempe, and D. M\"uller, {\it Estimates for maximal functions associated
to hypersurfaces in $\R^3$ and related problems of harmonic analysis}, Acta Math. {\bf 204} (2010), no. 2,
151--271.

\noindent [IM] I. Ikromov, D. M\"uller, {\it  Uniform estimates for the Fourier transform of surface-carried measures in
$\R^3$ and an application to Fourier restriction}, J. Fourier Anal. Appl. {\bf 17} (2011), no. 6, 1292-1332.

\noindent [Io]  A. Iosevich, {\it Maximal operators associated to families of flat curves in the plane},
Duke Math. J. {\bf 76} no. 2 (1994) 633-644. 

\noindent [IoSa] A. Iosevich, E. Sawyer, {\it Maximal averages over surfaces},  Adv. Math. {\bf 132} 
(1997), no. 1, 46--119.

\noindent [NSeW] A. Nagel, A. Seeger, and S. Wainger, {\it Averages over convex hypersurfaces},
Amer. J. Math. {\bf 115} (1993), no. 4, 903--927.

\noindent [PS] D. H. Phong and E. M. Stein, {\it The Newton polyhedron and oscillatory integral operators},
Acta Math. {\bf 179} (1997), no. 1, 105-152. 

\noindent [Se] A. Seeger, {\it Radon transforms and finite type conditions}, J. Amer. Math. Soc. {\bf 11} (1998), no. 4, 869-897. 

\noindent [So] C. Sogge, {\it Maximal operators associated to hypersurfaces with one nonvanishing princi pal 
curvature}  (English summary)  in {\it Fourier analysis and partial differential equations} (Miraflores de 
la Sierra, 1992),  317--323, Stud. Adv. Math., CRC, Boca Raton, FL, 1995. 

\noindent [SoS] C. Sogge and E. Stein, {\it Averages of functions over hypersurfaces in $\R^n$}, Invent. Math. {\bf 82} (1985), no. 3, 543--556.

\noindent [S1] E. Stein, {\it Maximal functions. I. Spherical means.} Proc. Nat. Acad. Sci. U.S.A. {\bf 73} (1976), no. 7, 2174--2175. 

\noindent [S2] E. Stein, {\it Harmonic analysis; real-variable methods, orthogonality, and oscillatory 
integrals}, Princeton Mathematics Series Vol. 43, Princeton University Press, Princeton, NJ, 1993.

\noindent [St] B. Street, {\it Sobolev spaces associated to singular and fractional Radon transforms}, Rev. Mat. Iberoam. {\bf 33} (2017), no. 2, 633-748.

\noindent [V] A. N. Varchenko, {\it Newton polyhedra and estimates of oscillatory integrals}, Functional 
Anal. Appl. {\bf 18} (1976), no. 3, 175-196.

\noindent [Z] E. Zimmermann, {\it On $L^p$-estimates for maximal averages over hypersurfaces not satisfying the transversality condition},
PhD thesis, Kiel University, 2014.

\
\

\noindent Department of Mathematics, Statistics, and Computer Science \hfill \break
\noindent University of Illinois at Chicago \hfill \break
\noindent 322 Science and Engineering Offices \hfill \break
\noindent 851 S. Morgan Street \hfill \break
\noindent Chicago, IL 60607-7045 \hfill \break
\noindent greenbla@uic.edu \hfill\break

\end{document}